\documentclass[final]{siamart171218}
\usepackage[todo]{style}
\usepackage{cleveref}
\usepackage{mathtools}

\headers{Consistent Second Moment Methods}{Olivier, Southworth, Warsa, and Park}
\title{{Consistent Second Moment Methods with Scalable Linear Solvers for Radiation Transport}\thanks{B.S.S., J.S.W., and H.K.P. were supported by the Laboratory Directed Research and Development program of Los Alamos National Laboratory under project number 20220174ER. S.O. was supported by the U.S. Department of Energy as a Nicholas C. Metropolis Fellow under the Laboratory Directed Research and Development program of the Los Alamos National Laboratory. Los Alamos National Laboratory report number LA-UR-24-23582. }}
\author{
    Samuel Olivier\thanks{Los Alamos National Laboratory, Computer, Computational, and Statistical Sciences Division, Los Alamos, New Mexico 87545 (\email{solivier@lanl.gov}).}
    \and
    Ben S. Southworth\thanks{Los Alamos National Laboratory, Theoretical Division, Los Alamos, New Mexico 87545.}
    \and
    James S. Warsa\footnotemark[2]
    \and
    HyeongKae Park\footnotemark[3]
}

\begin{document}
\maketitle

\begin{abstract}
Second Moment Methods (SMMs) are developed that are consistent with the Discontinuous Galerkin (DG) spatial discretization of the discrete ordinates (or \Sn) transport equations. 
The low-order (LO) diffusion system of equations is discretized with fully consistent \Pone, Local Discontinuous Galerkin (LDG), and Interior Penalty (IP) methods.
A discrete residual approach is used to derive SMM correction terms that make each of the LO systems consistent with the high-order (HO) discretization. 
We show that the consistent methods are more accurate and have better solution quality than independently discretized LO systems, that they preserve the diffusion limit, and that the LDG and IP consistent SMMs can be scalably solved in parallel on a challenging, multi-material benchmark problem. 
\end{abstract}


\section{Introduction}
The Second Moment Method (SMM) is an iterative scheme for solving high-dimensional kinetic equations that utilizes a low-order (LO) model to accelerate the calculation of slow-to-converge physics such as scattering in linear radiation transport and absorption-emission in thermal radiative transfer \cite{lewis_miller}. 
SMM is a member of a broad class of methods that includes the Variable Eddington Factor (VEF) \cite{mihalas}, Quasidiffusion (QD) \cite{goldin}, and high-order low-order (HOLO) \cite{doi:10.1080/00411450.2012.671224} methods. 
Such methods are particularly well suited for multiphysics calculations as the LO system can be directly coupled to the other physics in place of the high-order (HO) system, reducing reliance on expensive high-dimensional operations. 
Furthermore, this class of methods has been shown to be efficient across a wide range of problems including nuclear reactor modeling \cite{doi:10.13182/NSE11-81}, ocean modeling \cite{10.1016/j.procs.2015.05.477}, and gas dynamics \cite{CHACON201721}. 

A moment method is called ``algebraically consistent'' if the discrete LO solution matches the moments of the discrete HO solution to iterative solver tolerances independently from the mesh size.
By tying the LO solution to the moments of the HO solution, consistency guarantees that the LO system reproduces the solution quality of the HO solution and that both HO and LO solutions are conservative. 
The requirement of consistency often results in LO discretizations that can make the application of existing linear solver technology difficult or impractical \cite{warsa_mfem,dima_dfem}. 
In contrast, relaxing the consistency restraint leads to so-called independent methods \cite{two-level-independent-warsa} which allow efficient solvers to be applied to the LO system \cite{dgvef_olivier,me,rtvef_olivier,olivier_smm}, leading to greater overall computational efficiency. 
However, in these schemes, only the LO solution is conservative requiring extra care when multiphysics coupling beyond that which can be provided by the LO system alone is needed. 
Furthermore, it was shown in \citet{dgvef_olivier} that the LO solution was slow to converge to the HO solution, requiring significant mesh refinement to match the solution quality produced by the HO scheme in isolation. 
The independent methods from \cite{dgvef_olivier,me,rtvef_olivier,olivier_smm} also produced suboptimally converging approximations for the first moment. 
Finally, consistency was recognized as important in a recently developed HOLO algorithm with implicit-explicit time integration \cite{imex-trt}. 
Our results support the importance of consistent LO discretizations, particularly for the first moment where independent LO systems can yield very poor solution quality despite converging to the HO system with mesh refinement (see Figs.~\ref{fig:line_J} and \ref{fig:indep_ho_v_lo}).

In this paper, we design discretizations of the SMM low-order system that can be scalably solved and are algebraically consistent with the Discontinuous Galerkin (DG) spatial discretization of the Discrete Ordinates (\Sn) transport equations. 
In particular, we derive consistent SMMs with LO systems equivalent to the fully consistent \Pone \cite{WWM}, Local Discontinuous Galerkin (LDG) \cite{10.1007/s10915-007-9130-3}, and Interior Penalty (IP) \cite{Arnold2002} discretizations of radiation diffusion. 
Our approach is to use the discrete residual procedure developed for HOLO methods \cite{doi:10.1080/00411450.2012.671224,doi:10.1080/00295639.2020.1769390} -- where the LO correction terms are formed as a discrete residual that force a given LO system to match the moments of HO -- to generate the correction source terms needed for consistency. 
Through this discrete residual, consistency with the HO discretization can be achieved for any LO discretization, allowing the flexibility to choose LO discretizations that have known, efficient linear solvers without sacrificing consistency.

Our motivation for this combination of schemes is that the DG \Sn transport discretization is commonly used in thermal radiative transfer calculations and the LDG and IP discretizations are well known to be effectively preconditioned by Algebraic Multigrid (AMG).
The fully consistent \Pone discretization is included as a comparison method as it is known to be remarkably effective when used in the context of Diffusion Synthetic Acceleration (DSA) but notoriously difficult to solve iteratively \cite{warsa_mfem}, leaving unscalable, sparse direct methods as the only viable alternative. 
In addition, each of the proposed LO discretizations produce the zeroth and first moments in the same spaces as the zeroth and first moments of the DG \Sn transport equation, avoiding the projections and interpolations needed to combine DG \Sn HO with finite volume \cite{doi:10.1080/00411450.2012.671224,doi:10.1080/00295639.2020.1769390} or vertex-centered \cite{WLA} LO discretizations. 
In this sense, the LDG and IP LO systems derived here are as close as possible to the fully consistent \Pone method while maintaining the solvability of the LO system. 

The paper proceeds as follows. 
In \Cref{sec:2mm} we describe the SMM algorithm for the continuous transport system and in \Cref{sec:transport_disc} define the DG \Sn transport discretization. 
We then derive the discrete moments of the DG \Sn transport equations and form the discrete residual for each of the \Pone, LDG, and IP LO systems in \Cref{sec:cons-2mm}. 
The resulting closures are analyzed and compared to the SMM closures of analogous independent SMMs from \cite{olivier_smm}. 
\Cref{sec:results} provides numerical results verifying the order of accuracy of the schemes, their iterative efficiency in the thick diffusion limit, and their weak scaling performance on a challenging, multi-material benchmark problem. 
Finally, \Cref{sec:conc} gives conclusions. 

\section{Second Moment Method Algorithm} \label{sec:2mm}
We consider the steady-state, mono-energetic, linear transport problem with isotropic scattering given by
\begin{subequations}
\begin{equation} \label{eq:transport}
	\Omegahat\cdot\nabla\psi + \sigma_t \psi = \frac{\sigma_s}{4\pi}\varphi + q \,, \quad \x \in \D \,,  
\end{equation}
\begin{equation}
	\psi(\x,\Omegahat) = \bar{\psi}(\x,\Omegahat) \,, \quad \x \in \partial \D \ \text{and} \ \Omegahat\cdot\n < 0 \,, 
\end{equation}
\end{subequations}
where $\D \subset \R^{\dim}$ is the domain of the problem with $\partial \D$ its boundary, $\x\in \D$ and $\Omegahat\in \mathbb{S}^2$ are the spatial and angular variables, respectively, $\psi(\x,\Omegahat)$ the angular flux, $\sigma_t(\x)$ and $\sigma_s(\x)$ the total and scattering macroscopic cross sections, respectively, $q(\x,\Omegahat)$ the fixed-source, and $\bar{\psi}(\x,\Omegahat)$ a specified inflow boundary function. 
We define the zeroth, first, and second moments of the angular flux as $\varphi = \int \psi \ud\Omega$, $\vec{J} = \int \Omegahat\,\psi \ud \Omega$, and $\P = \int \Omegahat\otimes\Omegahat\,\psi \ud \Omega$ and refer to them as the scalar flux, current, and pressure, respectively. 
This transport problem is often extended to time and energy dependence and coupled with more complicated physics models such as thermal radiative transfer and radiation-hydrodynamics. 
In SMM, the scattering term is determined by the scalar flux solution of a certain LO system which we now derive. 
The zeroth and first angular moments of \eqref{eq:transport} are 
	\begin{subequations}
	\begin{equation}
		\nabla\cdot\vec{J} + \sigma_a \varphi = Q_0 \,, 
	\end{equation}
	\begin{equation}
		\nabla\cdot \P + \sigma_t \vec{J} = \vec{Q}_1 \,,
	\end{equation}
	\end{subequations}
where $\sigma_a(\x) = \sigma_t(\x) - \sigma_s(\x)$ is the absorption macroscopic cross section and $Q_i = \int \Omegahat^i\, q \ud \Omega$ are the angular moments of the transport equation's fixed-source, $q$. 
Let $J_n^\pm = \int_{\Omegahat\cdot\n\gtrless 0} \Omegahat\cdot\n\, \psi \ud \Omega$ denote the partial currents. 
Suitable boundary conditions are found by manipulating partial currents according to 
	\begin{equation}
		\vec{J}\cdot\n = J_n^+ + J_n^- = (J_n^+ - J_n^-) + 2J_n^- = B(\psi) + 2 J_n^- \,, 
	\end{equation}
where $B(\psi) = \int |\Omegahat\cdot\n|\,\psi \ud \Omega$. 
The moments of the transport equation have more unknowns than equations. 
Thus, we introduce auxiliary equations, referred to as closures, that relate $\P$ and $B$ in terms of the scalar flux and current. 
SMM uses additive closures of the form 
	\begin{subequations}
	\begin{equation}
		\P = \frac{1}{3}\I \varphi + \T(\psi) \,, 
	\end{equation}
	\begin{equation}
		B = \frac{1}{2}\varphi + \beta(\psi) \,,
	\end{equation}
where $\I$ denotes the $\dim\times\dim$ identity tensor and 
	\end{subequations}
	\begin{subequations} \label{eq:smm_closures}
	\begin{equation} \label{eq:smm_closures_T}
		\T(\psi) = \int \paren{\Omegahat\otimes\Omegahat - \frac{1}{3}\I}\psi \ud \Omega \,,
	\end{equation}
	\begin{equation} \label{eq:smm_closures_beta}
		\beta(\psi) = \int \paren{|\Omegahat\cdot\n| - \frac{1}{2}}\psi \ud \Omega \,,
	\end{equation}
	\end{subequations}
are the SMM correction tensor and boundary factor, respectively. 
Observe that the SMM closures are simply algebraic manipulations of $\P$ and $B$, respectively. 
Thus, when the angular flux is known, the LO SMM system is an equivalent reformulation of the zeroth and first angular moments of the HO transport equation. 

An efficient iterative algorithm is found by alternating HO and LO solves: the left hand side of the transport equation is inverted against a fixed scattering source, the closures $\T$ and $\beta$ are computed, and the LO system is then solved to update the transport equation's scattering source. 
The closures are weak functions of the HO solution so even simple iterative schemes converge rapidly and robustly with a cost that is independent of the mean free path. 

\section{Transport Discretization} \label{sec:transport_disc}
The Discrete Ordinates (\Sn) method is used to collocate the direction variable at a set of discrete points from a quadrature rule on the unit sphere $\{w_d,\Omegahat_d\}_{d=1}^{N_\Omega}$ such that $\psi_d(\x) = \psi(\x,\Omegahat_d)$.
Each discrete angle is then discretized in space with linear Discontinuous Galerkin (DG) finite elements so that, for each direction $\Omegahat_d$, $\psi_d \in Y_1$ where $Y_1$ is the space of piecewise discontinuous linear polynomials. 
Let the domain $\D$ be split into finite elements, $K$, such that $\D = \cup K$. 
The spatial discretization is derived by multiplying the transport equation for each $\Omegahat_d$ by an arbitrary test function $u \in Y_1$ and integrating over each element in the mesh. 
On interior interfaces, upwinding is used such that 
	\begin{equation}
		\Omegahat_d\cdot\n\, \widehat{\psi}_d = \begin{cases}
			\Omegahat_d\cdot\n\, \psi_{d,1} \,, & \Omegahat_d\cdot\n > 0 \\
			\Omegahat_d\cdot\n\, \psi_{d,2} \,, & \Omegahat_d\cdot\n < 0 
		\end{cases} \,,
	\end{equation}
where $\psi_{d,i} = \psi_d|_{K_i}$ and the indices ``1'' and ``2'' correspond to two arbitrary neighboring elements in the mesh following the convention that the normal vector, $\n$, points from $K_1$ to $K_2$. 
Note that the upwind flux can be equivalently written using the switch functions
	\begin{equation}
		\Omegahat_d\cdot\n\, \widehat{\psi}_d = \frac{1}{2}\paren{\Omegahat_d\cdot\n + |\Omegahat_d\cdot\n|}\psi_{d,1} + \frac{1}{2}\paren{\Omegahat_d\cdot\n - |\Omegahat_d\cdot\n|}\psi_{d,2} \,. 
	\end{equation}
Combining like terms, this is equivalent to 
	\begin{equation}\label{eq:upwind}
		\Omegahat_d\cdot\n\, \widehat{\psi}_d = \frac{\Omegahat\cdot\n}{2}(\psi_{d,1} + \psi_{d,2}) + \frac{|\Omegahat\cdot\n|}{2}(\psi_{d,1} - \psi_{d,2}) = \Omegahat\cdot\n\avg{\psi_d} + \frac{|\Omegahat\cdot\n|}{2}\jump{\psi_d} 
	\end{equation}
where, for $u \in Y_1$, 
	\begin{equation}
		\jump{u} = u_1 - u_2 \,, \quad \avg{u} = \frac{1}{2}(u_1 + u_2) \,, 
	\end{equation}
are the jump and average, respectively, with an analogous definition for vector-valued arguments. 
Note that while $|\Omegahat\cdot\n|$ is continuous, it has a discontinuous derivative at $\Omegahat\cdot\n = 0$ and thus angular moments of the upwind numerical flux will be difficult to compute accurately with numerical quadrature. 
On the boundary of the domain, the numerical flux is 
	\begin{equation}
		\Omegahat_d\cdot\n\, \widehat{\psi}_d = \begin{cases}
			\Omegahat_d\cdot\n\, \psi_{d} \,, & \Omegahat_d\cdot\n > 0 \\ 
			\Omegahat_d\cdot\n\, \bar{\psi}(\x,\Omegahat_d) \,, & \Omegahat_d\cdot\n < 0 
		\end{cases} \,. 
	\end{equation}
That is, for an element on the boundary, the inflow is specified by the boundary function, $\bar{\psi}$, and the outflow is specified by the numerical solution. 

With these definitions, the transport discretization is: for each $\Omegahat_d$, find $\psi_d \in Y_1$ such that 
\begin{multline} \label{eq:dgsn_transport}
	\int_{\Gamma_0} \Omegahat\cdot\n\, \jump{u}\avg{\psi_d} \ud s + \frac{1}{2}\int_{\Gamma_0} |\Omegahat\cdot\n|\, \jump{u}\jump{\psi_d} \ud s + \int_{\Gamma_{b,d}^+} \Omegahat\cdot\n\, u \psi \ud s - \int \Omegahat\cdot\nablah u\, \psi_d \ud \x \\+ \int \sigma_t\, u\psi_d \ud \x = \frac{1}{4\pi}\int \sigma_s\, u\varphi \ud \x + \int u\, q \ud \x - \int_{\Gamma_{b,d}^-} \Omegahat\cdot\n\, u \bar{\psi} \ud s \,, \quad \forall u \in Y_1 \,, 
\end{multline}
where $\Gamma_0$ is the set of unique interior faces in the mesh, $\Gamma_{b,d}^\pm$ are the outflow/inflow portions of the boundary of the domain, and $\nablah$ the element-local gradient. 
In the SMM algorithm, the scalar flux in the scattering term is determined by the LO system and thus each HO solve involves inverting the left hand side of \eqref{eq:dgsn_transport} only. 

\section{Derivation of Consistent SMMs}\label{sec:cons-2mm}
SMMs consistent with DG \Sn transport are derived by introducing an SMM correction term as a residual between the moments of the DG \Sn transport equations and a chosen LO discretization. 
Upon iterative convergence, terms present in the LO system but not in the HO system are canceled by this residual, making the solution of the LO system equivalent to the angular moments of the discrete transport equation. 
This process allows us to choose discretizations of the LO system for which scalable and efficient iterative solution methods are available while still maintaining consistency between the HO and LO equations. 

In this section, we derive the zeroth and first angular moments of the discrete transport equation and the correction source terms that make the fully consistent \Pone, Local Discontinuous Galerkin (LDG), and Interior Penalty (IP) discretizations of radiation diffusion consistent with the moments of DG \Sn transport. 
Additionally, the correction terms for the consistent IP method will be analyzed and compared to the independent variant of IP in \cite{olivier_smm}. 

\subsection{Angular Moments of DG \Sn Transport}
To notationally separate the HO and LO variables, we use the ``HO'' subscript to denote quantities computed from the HO solution, $\psi_d$. 
That is, the zeroth, first, and second moments of the HO solution are denoted
	\begin{equation}
		\phiho = \sum_d w_d\,\psi_d \,, \quad \Jho = \sum_d w_d\,\Omegahat_d\,\psi_d \,, \quad \Pho = \sum_d w_d\,\Omegahat_d\otimes\Omegahat_d\,\psi_d \,. 
	\end{equation}
Since each $\psi_d \in Y_1$, we have that $\phiho$ and each component of $\Jho$ and $\Pho$ are all members of $Y_1$. 
We write $[Y_1]^{\dim}$ to denote the $\dim$-fold tensor product of the space $Y_1$ such that $\Jho \in [Y_1]^{\dim}$ has each component in $Y_1$. 
We focus primarily on the upwind numerical fluxes in the discretized transport equation because the remaining terms are derived by simply replacing the discrete moments of $\psi_d$ with $\phiho$, $\Jho$, and $\Pho$, where appropriate. 
The zeroth and first moments of the upwind numerical flux, denoted by $\widehat{\vec{J}}_{\smallHO}\cdot\n$ and $\widehat{\P}_{\smallHO}\n$, respectively, require half-range integrations because the upwind flux cannot be integrated exactly by \Sn quadrature. 
From \eqref{eq:upwind}, the zeroth moment of the upwind numerical flux is
	\begin{equation} \label{eq:dgsn_Jflux}
	\begin{aligned}
		\widehat{\vec{J}}_{\smallHO}\cdot\n &= \sum_d w_d\,\Omegahat_d\cdot\n\, \widehat{\psi}_d \\
		&= \sum_d w_d \paren{\Omegahat_d\cdot\n \avg{\psi_d} + \frac{1}{2}|\Omegahat\cdot\n|\jump{\psi_d}} \\
		&= \avg{\sum_d w_d\,\Omegahat_d\cdot\n\, \psi_d} + \frac{1}{2}\jump{\sum_d w_d\,|\Omegahat_d\cdot\n|\,\psi_d} \\
		&= \avg{\Jho\cdot\n} + \frac{1}{2}\jump{\Jhoplus - \Jhominus} \,,
	\end{aligned}
	\end{equation}
where
	\begin{equation} \label{eq:Jhopm}
		\Jhopm = \sum_{\Omegahat_d\cdot\n\gtrless 0} w_d\,\Omegahat_d\cdot\n\,\psi_d 
	\end{equation}
are the half-range or ``partial'' currents. 
Applying the definitions of the jump and average and dropping the ``HO'' subscript for brevity, the zeroth moment of the upwind numerical flux is simply 
	\begin{equation} \label{eq:upwind_partial_current}
	\begin{aligned}
		\widehat{\vec{J}}\cdot\n &= \avg{\vec{J}\cdot\n} + \frac{1}{2}\jump{J_n^+ - J_n^-} \\
		&= \frac{1}{2}\paren{\vec{J}_1\cdot\n + \vec{J}_2\cdot\n} + \frac{1}{2}\paren{J_{n,1}^+ - J_{n,1}^- - (J_{n,2}^+ - J_{n,2}^-)} \\
		&= \frac{1}{2}\paren{J_{n,1}^+ + J_{n,1}^- + J_{n,2}^+ + J_{n,2}^-} + \frac{1}{2}\paren{J_{n,1}^+ - J_{n,1}^- - (J_{n,2}^+ - J_{n,2}^-)} \\
		&= J_{n,1}^+ + J_{n,2}^- \,,
	\end{aligned}
	\end{equation}
where $J_{n,i}^\pm = J_{n}^\pm|_{K_i}$ are the half-range currents computed using information in element $K_i$. 
In other words, the zeroth moment of the upwind numerical flux is computed from the sum of two partial currents, one in which information from element ``1'' is used to compute the particle flow from $K_1$ to $K_2$ and the other in which information from element ``2'' computes the flow in the opposite direction from $K_2$ to $K_1$. 
This viewpoint lends itself to a simpler implementation while \eqref{eq:dgsn_Jflux} is useful in analyzing the closures that follow in subsequent sections. 
On the boundary, the zeroth moment of the numerical flux is 
	\begin{equation} \label{eq:dgsn_Jflux_bdr}
		\widehat{\vec{J}}_{\smallHO}\cdot\n = \sum_d w_d\, \Omegahat_d\cdot\n\, \widehat{\psi}_d = \Jhoplus + J_\text{in} 
	\end{equation}
where 
	\begin{equation}
		J_\text{in} = \sum_{\Omegahat_d\cdot\n<0} w_d\,\Omegahat_d\cdot\n\,\bar{\psi}(\x,\Omegahat_d)
	\end{equation}
is the inflow partial current computed from the inflow boundary function, $\bar{\psi}$. 

Similarly, the first moment of the upwind numerical flux is
	\begin{equation} \label{eq:dgsn_Pflux}
		\widehat{\P}_{\smallHO}\n = \sum_d w_d\, \Omegahat_d(\Omegahat_d\cdot\n)\, \widehat{\psi}_d = \avg{\Pho\n} + \frac{1}{2}\jump{\Phoplus - \Phominus} \,,
	\end{equation}
where 
	\begin{equation}
		\Phopm = \sum_{\Omegahat_d\cdot\n\gtrless 0} w_d\,\Omegahat_d (\Omegahat_d\cdot\n) \, \psi_d 
	\end{equation}
are the half-range pressures. 
Again, applying the definitions of the jump and average yields $\widehat{\P}_{\smallHO}\n = \vec{P}_{\text{HO},n,1}^+ + \vec{P}_{\text{HO},n,2}^-$ so that the first moment of the upwind numerical flux is computed from the half range pressures. 
On the boundary, 
	\begin{equation} \label{eq:dgsn_Pflux_bdr}
		\widehat{\P}_{\smallHO}\n = \sum_d w_d\, \Omegahat_d(\Omegahat_d\cdot\n)\, \widehat{\psi}_d = \Phoplus + \vec{P}_\text{in} 
	\end{equation}
where  
	\begin{equation}
		\vec{P}_\text{in} = \sum_{\Omegahat_d\cdot\n<0} w_d\,\Omegahat_d(\Omegahat_d\cdot\n)\,\bar{\psi}(\x,\Omegahat_d) 
	\end{equation}
is the half-range pressure computed from the inflow transport boundary condition. 

The discrete zeroth moment is derived by operating on the DG \Sn transport equation from \eqref{eq:dgsn_transport} by the quadrature sum $\sum_d w_d\paren{\cdot}$. 
Each component of the discrete first moment is similarly derived by operating on \eqref{eq:dgsn_transport} with the quadrature sum $\sum_d w_d\,\e_i \cdot \Omegahat_d$ where $\e_i$ are the coordinate-axis unit vectors and $1 \leq i \leq \dim$. 
We use the fact that 
	\begin{equation}
		\{ \vec{v} : \vec{v} = \e_i u \quad \forall u \in Y_1 \ \mathrm{and} \ 1 \leq i \leq \dim \} 
	\end{equation}
forms a basis for the space $[Y_1]^{\dim}$ to write all the components of the first moment equation as a single equation with a vector-valued test function $\vec{v} = \e_i u$ where $u \in Y_1$. 
The zeroth and first moments of the DG \Sn transport discretization are 
	\begin{subequations} \label{eq:dgsn}
	\begin{multline} \label{eq:dgsn_zeroth}
		\int_{\Gamma_0} \jump{u}\paren{\avg{\Jho\cdot\n} + \frac{1}{2}\jump{\Jhoplus - \Jhominus}} \ud s + \int_{\Gamma_b} u\,\Jhoplus \ud s \\- \int \nablah u \cdot \Jho \ud \x + \int \sigma_a\, u\phiho \ud \x = \int u\, Q_0 \ud \x - \int_{\Gamma_b} u\, J_\text{in} \ud s \,,
	\end{multline}
	\begin{multline} \label{eq:dgsn_first}
		\int_{\Gamma_0} \jump{\vec{v}}\cdot\paren{\avg{\Pho\n} + \frac{1}{2}\jump{\Phoplus - \Phominus}} \ud s + \int_{\Gamma_b} \vec{v}\cdot\Phoplus\ud s \\ - \int \nablah \vec{v} : \Pho \ud \x + \int \sigma_t\, \vec{v}\cdot\Jho \ud \x = \int \vec{v}\cdot\vec{Q}_1 - \int \vec{v}\cdot\vec{P}_\text{in} \ud s \,,
	\end{multline}
	\end{subequations}
where $u \in Y_1$ and $\vec{v} \in [Y_1]^{\dim}$ are arbitrary and the angular moments of the source $q$ are computed with \Sn  quadrature. 
In the following, \eqref{eq:dgsn} will be algebraically manipulated so that the left hand side is equivalent to a desired LO discretization of radiation diffusion. 
Manipulations applied on the left hand side are balanced on the right hand side, producing correction terms that make the LO system equivalent to \eqref{eq:dgsn}. 

\subsection{\Pone}
The \Pone diffusion system is derived by restricting the transport solution to be at most linearly anisotropic. 
That is, the half-range current and half and full-range pressures are closed under the assumption that $\psi_d = \frac{1}{4\pi}\!\paren{\varphi + 3\Omegahat_d\cdot\vec{J}}$. 
Under this ansatz, the pressure is
	\begin{equation} \label{eq:diffusion}
		\P = \frac{1}{3} \I \varphi \,, 
	\end{equation}
and the half-range currents are
	\begin{equation}
		J_n^\pm = \frac{\vec{J}\cdot\n}{2} \pm \frac{\alpha}{2} \varphi \,. 
	\end{equation}
Note that using the quadrature-dependent term $\alpha = \sum_d w_d\,|\Omegahat_d\cdot\n|$ such that $\alpha \rightarrow 1/2$ as $N_\Omega \rightarrow \infty$ can be important for achieving optimal convergence for independent SMMs \cite{olivier_smm}. 
Using the upwinded partial current definition of $\widehat{\vec{J}}_{\smallHO}\cdot\n$ in \eqref{eq:upwind_partial_current} and the linearly anisotropic partial currents, the zeroth moment of the upwind numerical flux is 
	\begin{equation} \label{eq:P1_J}
	\begin{aligned} 
		\widehat{\vec{J}}\cdot\n &= J_{n,1}^+ + J_{n,2}^- \\
		&= \frac{\vec{J}_1\cdot\n}{2} + \frac{\alpha}{2}\varphi_1 + \frac{\vec{J}_2\cdot\n}{2} - \frac{\alpha}{2} \varphi_2 \\
		&= \avg{\vec{J}\cdot\n} + \frac{\alpha}{2}\jump{\varphi} \,. 
	\end{aligned}
	\end{equation}
The first moment is closed using \eqref{eq:diffusion} and the partial currents are used to find an expression for the scalar flux. 
That is, 
	\begin{equation}
		J_n^+ - J_n^- = \alpha \varphi \Rightarrow \varphi = \frac{1}{\alpha}(J_n^+ - J_n^-) \,. 
	\end{equation}
This expression is combined with \eqref{eq:diffusion} and upwinded so that the first moment of the upwind numerical flux is:
	\begin{equation} \label{eq:P1_P}
	\begin{aligned} 
		\widehat{\P}\n &= \frac{\n}{3\alpha}\paren{J_{n,1}^+ - J_{n,2}^-} \\
		&= \frac{\n}{3\alpha}\bracket{\frac{\vec{J}_1\cdot\n}{2} + \frac{\alpha}{2}\varphi_1 - \paren{\frac{\vec{J}_2\cdot\n}{2} - \frac{\alpha}{2}\varphi_2}} \\
		&= \frac{\n}{6\alpha}\jump{\vec{J}\cdot\n} + \frac{\n}{3}\avg{\varphi} \,. 
	\end{aligned}
	\end{equation}
Because $\widehat{\P}\n$ depends on $\vec{J}\cdot\n$, this numerical flux couples the current degrees of freedom across interior faces. 
For boundary faces, information from the ``2'' element is computed from the half range moment of the inflow source function so that 
	\begin{subequations} \label{eq:lobc}
	\begin{equation} \label{eq:lobc_J}
		\widehat{\vec{J}}\cdot\n = \frac{\vec{J}\cdot\n}{2} + \frac{\alpha}{2}\varphi + J_\text{in} \,, 
	\end{equation}
	\begin{equation} \label{eq:lobc_P}
		\widehat{\P}\n = \frac{\n}{6\alpha}\vec{J}\cdot\n + \frac{\n}{6}\varphi + \vec{P}_\text{in} \,. 
	\end{equation}
	\end{subequations}
The moments of DG \Sn transport are manipulated to make the left hand side of \eqref{eq:dgsn} equivalent to \Pone. 
To that end, we seek to define the numerical fluxes and pressure such that, upon iterative convergence, the HO and LO variables cancel leaving only the numerical fluxes and pressure associated with the HO moments. 
First, consider $\widehat{\vec{J}}_{\smallHO}\cdot\n$ from \eqref{eq:dgsn_Jflux}. 
We add the \Pone numerical flux in \eqref{eq:P1_J} to $\widehat{\vec{J}}_{\smallHO}\cdot\n$ and subtract the \Pone numerical flux computed from the HO solution. 
The resulting \Pone SMM numerical flux is 
	\begin{equation} \label{eq:P1_Jhat}
	\begin{aligned}
		\widehat{\vec{J}}_\mathrm{P1}\cdot\n &= \widehat{\vec{J}}_{\smallHO}\cdot\n + \avg{\vec{J}\cdot\n} + \frac{\alpha}{2}\jump{\varphi} - \avg{\Jho\cdot\n} - \frac{\alpha}{2}\jump{\phiho} \\
		&= \avg{\vec{J}\cdot\n} + \frac{\alpha}{2}\jump{\varphi} + \frac{1}{2}\jump{\Jhoplus - \Jhominus - \alpha \phiho} \\
		&= \avg{\vec{J}\cdot\n} + \frac{\alpha}{2}\jump{\varphi} + \frac{1}{2}\jump{\beta} \,,
	\end{aligned}
	\end{equation}
where $\beta = \sum_d w_d\,(|\Omegahat_d\cdot\n| - \alpha)\,\psi_d$ is the discrete version of the SMM boundary correction factor from \eqref{eq:smm_closures_beta}. 
For boundary faces, an analogous process applied to the boundary numerical flux in \eqref{eq:dgsn_Jflux_bdr}, using the \Pone boundary flux from \eqref{eq:lobc_J}, yields 
	\begin{equation} \label{eq:P1_Jhat_bdr}
	\begin{aligned}
		\widehat{\vec{J}}_\mathrm{P1}\cdot\n &= \Jhoplus + J_\text{in} + \frac{\vec{J}\cdot\n}{2} + \frac{\alpha}{2}\varphi - \frac{\Jho}{2} - \frac{\alpha}{2} \phiho \\
		&= \frac{\vec{J}\cdot\n}{2} + \frac{\alpha}{2}\varphi + \frac{1}{2}\paren{\Jhoplus - \Jhominus - \alpha \phiho} + J_\text{in}\\
		&= \frac{\vec{J}\cdot\n}{2} + \frac{\alpha}{2}\varphi + \frac{\beta}{2} + J_\text{in} \,. 
	\end{aligned}
	\end{equation}
Repeating for the first moment of the numerical flux, now using \eqref{eq:P1_P} and \eqref{eq:dgsn_Pflux}, gives 
	\begin{equation}  \label{eq:P1_Phat}
	\begin{aligned}
		\widehat{\P}_\mathrm{P1}\n &= \widehat{\P}_{\smallHO}\n + \frac{\n}{6\alpha}\jump{\vec{J}\cdot\n} + \frac{\n}{3}\avg{\varphi} - \frac{\n}{6\alpha}\jump{\Jho\cdot\n} - \frac{\n}{3}\avg{\phiho} \\
		&= \frac{\n}{6\alpha}\jump{\vec{J}\cdot\n} + \frac{\n}{3}\avg{\varphi} + \avg{\Pho\n - \frac{\n}{3}\phiho} + \frac{1}{2}\jump{\Phoplus - \Phominus - \frac{\n}{3\alpha}\Jho\cdot\n} \\
		&= \frac{\n}{6\alpha}\jump{\vec{J}\cdot\n} + \frac{\n}{3}\avg{\varphi} + \avg{\T\n} + \frac{1}{2}\jump{\Phoplus - \Phominus - \frac{\n}{3\alpha}\Jho\cdot\n} 
	\end{aligned}
	\end{equation}
where $\T = \sum_d w_d\,(\Omegahat_d\otimes\Omegahat_d - 1/3\I)\, \psi_d$ is the discrete version of the SMM correction tensor from \eqref{eq:smm_closures_T}. 
On the boundary, \eqref{eq:dgsn_Pflux_bdr} is manipulated with \eqref{eq:lobc_P} to yield
	\begin{equation} \label{eq:P1_Phat_bdr}
	\begin{aligned}
		\widehat{\P}_\mathrm{P1}\n &= \Phoplus + \vec{P}_\text{in} + \frac{\n}{6\alpha}\vec{J}\cdot\n+ \frac{\n}{6}\varphi - \frac{\n}{6\alpha}\Jho\cdot\n - \frac{\n}{6}\phiho \,. 
	\end{aligned}
	\end{equation}
Finally, the full range pressure is closed with
	\begin{equation} \label{eq:SMM_P}
		\P_\text{SMM} = \Pho + \frac{1}{3}\I \varphi - \frac{1}{3}\I \phiho = \frac{1}{3}\I \varphi + \T \,. 
	\end{equation}
Observe that, if the HO and LO scalar flux and current were equal, the HO and LO variables in \eqref{eq:P1_Jhat}, \eqref{eq:P1_Jhat_bdr}, \eqref{eq:P1_Phat}, \eqref{eq:P1_Phat_bdr}, and \eqref{eq:SMM_P} would cancel leaving only the corresponding HO numerical flux or pressure, as desired.  

Substituting the \Pone interior and boundary numerical fluxes and the SMM pressure for their HO counterparts in \eqref{eq:dgsn} and moving all HO terms to the right hand side yields the LO \Pone SMM problem: find $(\varphi,\vec{J}) \in Y_1 \times [Y_1]^{\dim}$ such that
\begin{subequations} \label{eq:p1-weak}
\begin{multline} \label{eq:p1-m0}
	\int_{\Gamma_0} \jump{u} \avg{\vec{J}\cdot\n} \ud s + \frac{\alpha}{2}\int_{\Gamma_0} \jump{u}\jump{\varphi} \ud s + \frac{1}{2}\int_{\Gamma_b} u\, \vec{J}\cdot\n\ud s + \frac{\alpha}{2}\int_{\Gamma_b} u \varphi \ud s \\- \int \nablah u \cdot \vec{J} \ud \x + \int \sigma_a\, u \varphi \ud \x = \int u\, Q_0 \ud \x - \int_{\Gamma_b} u\, J_\text{in} \ud s + \mathcal{R}_\mathrm{P1}^0(u) \,, \quad \forall u \in Y_1 \,, 
\end{multline}
\begin{multline} \label{eq:p1-m1}
	\frac{1}{6\alpha}\int_{\Gamma_0} \jump{\vec{v}\cdot\n}\jump{\vec{J}\cdot\n} \ud s + \frac{1}{3}\int_{\Gamma_0} \jump{\vec{v}\cdot\n}\avg{\varphi} \ud s + \frac{1}{6\alpha}\int_{\Gamma_b} (\vec{v}\cdot\n)(\vec{J}\cdot\n) \ud s + \frac{1}{6}\int_{\Gamma_b} \vec{v}\cdot\n\, \varphi \ud s \\ - \frac{1}{3}\int \nablah\cdot\vec{v}\, \varphi \ud \x + \int \sigma_t\, \vec{v}\cdot\vec{J} \ud \x = \int \vec{v}\cdot\vec{Q}_1 \ud \x - \int_{\Gamma_b} \vec{v}\cdot\vec{P}_\text{in} \ud s + \mathcal{R}_\mathrm{P1}^1(\vec{v}) \,, \quad \forall \vec{v} \in [Y_1]^{\dim} \,, 
\end{multline}
\end{subequations}
where 
\begin{subequations}
\begin{equation}
	\mathcal{R}_\mathrm{P1}^0(u) = -\frac{1}{2}\int_{\Gamma_0} \jump{u}\jump{\beta} \ud s - \frac{1}{2}\int_{\Gamma_b} u\, \beta \ud s \,, 
\end{equation}
\begin{multline}
	\mathcal{R}_\mathrm{P1}^1(\vec{v}) = -\int_{\Gamma_0} \jump{\vec{v}} \cdot \avg{\T\n} \ud s - \frac{1}{2}\int_{\Gamma_0} \jump{\vec{v}}\cdot\jump{\Phoplus - \Phominus - \frac{\n}{3\alpha} \Jho\cdot\n} \ud s \\
	- \int_{\Gamma_b} \vec{v}\cdot\paren{\Phoplus - \frac{\n}{6\alpha}\Jho\cdot\n - \frac{\n}{6}\phiho} \ud s + \int \nablah \vec{v} : \T \ud \x \,, 
\end{multline} 
\end{subequations}
are \Pone correction source terms for the zeroth and first moments, respectively. 
Note that this system can be symmetrized by multiplying the first moment equation by $-3$. 
To show this, the local divergence is integrated by parts on each element:
	\begin{equation}
	\begin{aligned}
	 	\int \nablah\cdot\vec{v}\,\varphi \ud \x &= \int_{\Gamma_0} \jump{\vec{v}\cdot\n\,\varphi} \ud s + \int_{\Gamma_b} \vec{v}\cdot\n\,\varphi \ud s - \int \vec{v}\cdot\nablah \varphi \ud \x \\
	 	&= \int_{\Gamma_0} \jump{\vec{v}\cdot\n}\avg{\varphi} \ud s + \int_{\Gamma_0} \avg{\vec{v}\cdot\n}\jump{\varphi} \ud s + \int_{\Gamma_b} \vec{v}\cdot\n\, \varphi \ud s - \int \vec{v}\cdot\nablah\varphi \ud \x 
	\end{aligned}
	\end{equation} 
so that 
	\begin{multline}
		\frac{1}{3}\int_{\Gamma_0} \jump{\vec{v}\cdot\n}\avg{\varphi} \ud s + \frac{1}{6}\int_{\Gamma_b} \vec{v}\cdot\n\, \varphi \ud s - \frac{1}{3}\int \nablah\cdot\vec{v}\, \varphi \ud \x \\= -\frac{1}{3}\int_{\Gamma_0} \avg{\vec{v}\cdot\n}\jump{\varphi} \ud s - \frac{1}{6}\int_{\Gamma_b} \vec{v}\cdot\n\, \varphi \ud s + \frac{1}{3}\int \vec{v}\cdot\nablah \varphi \ud \x \,. 
	\end{multline}
The terms on the left hand side appear in the first moment equation \eqref{eq:p1-m1} and are then symmetric up to a factor of $-1/3$ with 
	\begin{equation}
		\int_{\Gamma_0} \jump{u} \avg{\vec{J}\cdot\n} \ud s + \frac{1}{2}\int_{\Gamma_b} u\,\vec{J}\cdot\n \ud s - \int \nablah u \cdot\vec{J} \ud \x 
	\end{equation}
from the zeroth moment. 
The remaining terms are symmetric. 

\subsection{Consistent Low-Order Systems with Scalable Solvers}
As mentioned, the \Pone LO system rapidly accelerates the solution of the HO system, but is difficult to solve iteratively. 
We now design LO systems so that the left hand side matches existing discretizations for elliptic problems to which efficient and scalable solvers can be applied, with right hand side source terms that result in consistent methods. 

After discretizing the LO system, the assembled matrices take a block form
\begin{equation}\label{eq:block}
A_\text{LO} = \begin{bmatrix} M_{J} & -1/3D^T \\ D & M_{\varphi} \end{bmatrix} \,. 
\end{equation}
We consider LO systems based on the LDG \cite{10.1007/s10915-007-9130-3} and IP \cite{Arnold2002} methods. 
Both schemes allow for an efficient solution procedure by severing the coupling between the normal component of the current across interior mesh interfaces. 
That is, the numerical flux for the first moment equation is chosen to be independent of $\vec{J}\cdot\n$. 
This makes $M_J$ block-diagonal by element, allowing the current to be eliminated directly to form an exact block LDU decomposition, with Schur complement in scalar flux, $S_{\varphi}= M_{\varphi} + 1/3D M_J^{-1}D^T$.
Note that $S_{\varphi}$ has the form of a discrete diffusion operator and can be scalably preconditioned with AMG. 
Once the scalar flux is known, the current solution can be easily computed with element-local back substitution.

\begin{figure}
\centering
\begin{subfigure}{.32\textwidth}
	\centering
	\includegraphics[width=\textwidth]{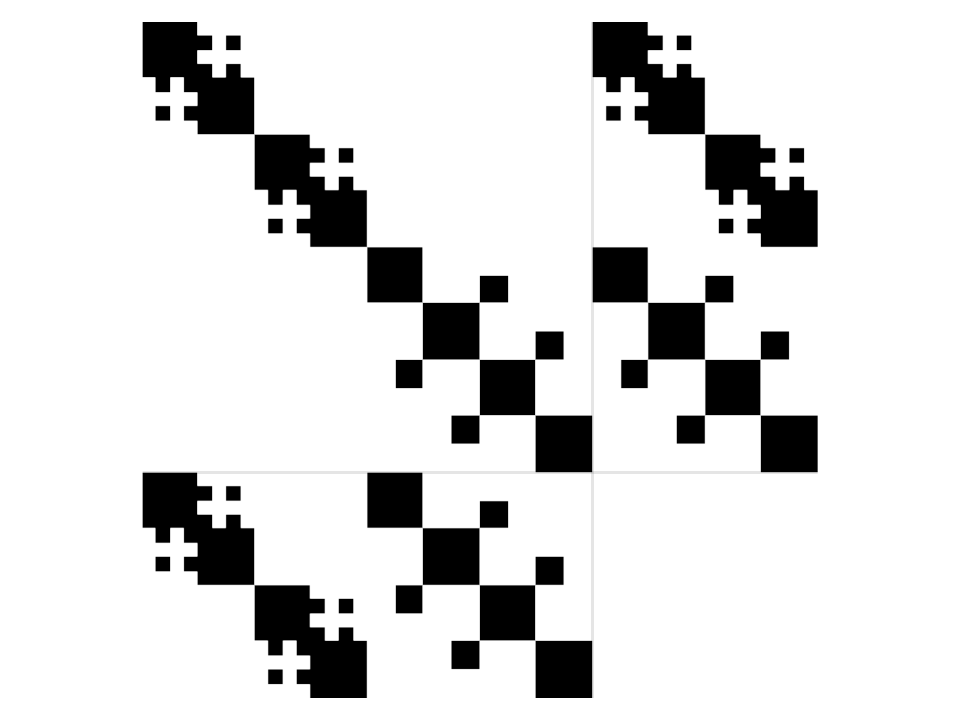}
	\caption{}
\end{subfigure}
\begin{subfigure}{.32\textwidth}
	\centering
	\includegraphics[width=\textwidth]{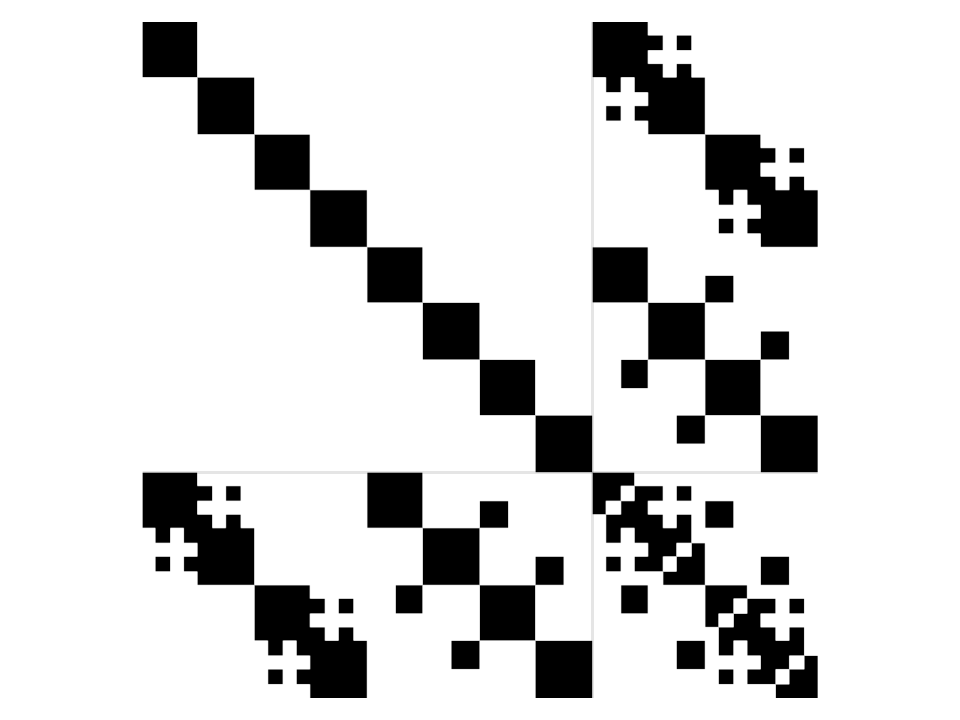}
	\caption{}
\end{subfigure}
\begin{subfigure}{.32\textwidth}
	\centering
	\includegraphics[width=\textwidth]{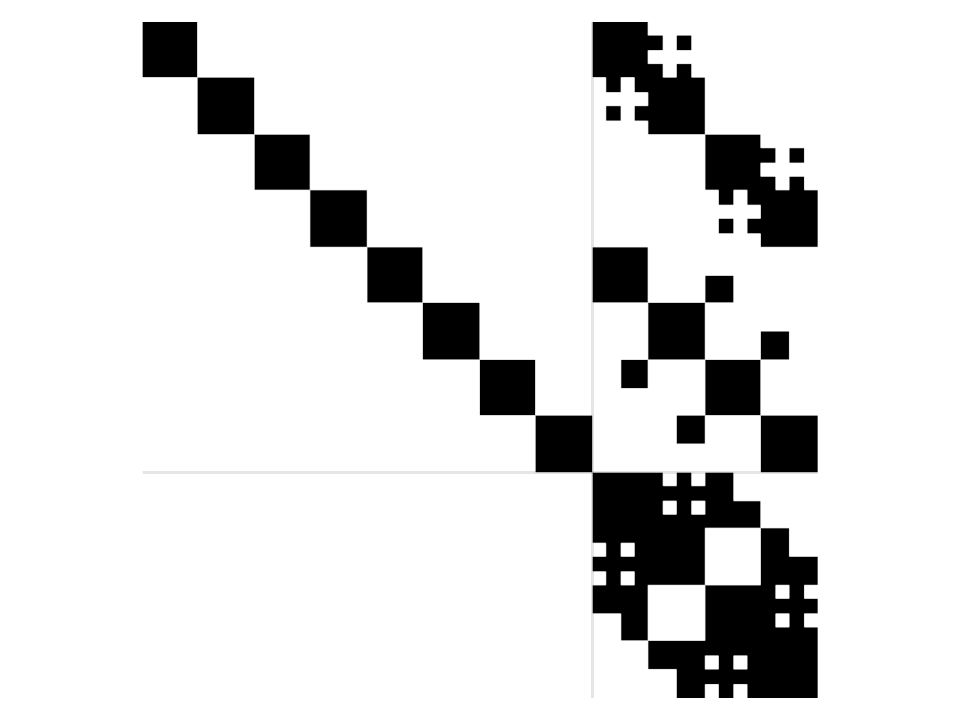}
	\caption{}
\end{subfigure}
\caption{Sparsity plots depicting DG discretizations of the Poisson equation in first-order form. The degrees of freedom are ordered such that the first block row corresponds to the current and the second to the scalar flux, as in \eqref{eq:block}. Plot (a) shows the sparsity pattern for a \Pone-like numerical flux. Observe that the top diagonal block is coupled to its neighbors, preventing the efficient elimination of the current. Plot (b) shows the sparsity pattern associated with the approach taken by the LDG and IP discretizations where the coupling in the normal component of the current is severed and replaced with alternate stabilization terms, allowing the current to be efficiently eliminated on each element. Plot (c) shows the sparsity pattern that results after eliminating the current with block Gaussian elimination. }
\label{fig:unified_spy}
\end{figure}

Due to the localization of the current, stabilization must be included to ensure the resulting discretizations are stable with respect to $h$.
LDG and IP differ only in their choice of stabilization. 
Figure \ref{fig:unified_spy} shows sparsity plots of the block IP system \eqref{eq:block} as the current is localized and stabilization is added. 
The stabilization terms are artifacts of the choice of the LO discretization.
Maintaining consistency requires modified SMM correction source terms that will impact the convergence rate of the SMM algorithm. 

\subsubsection{Local Discontinuous Galerkin}
We first consider the LDG method which is notable for its lack of problem-dependent parameters and avoidance of a penalty parameter that scales inversely with the mesh size. 
LDG is also naturally formulated and typically implemented as the Schur complement in the scalar flux of the block $2\times 2$ system \eqref{eq:block}. 
The LDG numerical fluxes are based on an arbitrary upwinding of the current that is balanced by an opposing choice for the scalar flux. 
Following \citet{10.1007/s10915-007-9130-3}, the LDG numerical fluxes are
	\begin{subequations}
	\begin{equation} \label{eq:ldg_J}
		\widehat{\vec{J}}\cdot\n = \avg{\vec{J}\cdot\n} + \frac{s}{2}\jump{\vec{J}\cdot\n} + \kappa \jump{\varphi} \,,
	\end{equation}
	\begin{equation} \label{eq:ldg_P}
		\widehat{\P}\n = \frac{\n}{3}\!\paren{\avg{\varphi} - \frac{s}{2}\jump{\varphi}} \,. 
	\end{equation}
	\end{subequations}
When $s=+1$, $\avg{\vec{J}\cdot\n} + s/2\jump{\vec{J}\cdot\n} = \vec{J}_1\cdot\n$ with $\vec{J}_2\cdot\n$ chosen when $s=-1$. 
The scalar flux is chosen oppositely, with $\varphi_2$ and $\varphi_1$ selected when $s = +1$ and $s=-1$, respectively. 
The upwinding parameter $s$ is chosen according to 
	\begin{equation}
		s = \begin{cases}
			+1 \,, & \vec{w}\cdot\n > 0 \\ 
			-1 \,, & \vec{w} \cdot\n < 0 
		\end{cases}
	\end{equation}
for some arbitrary, non-zero vector $\vec{w} \in \R^{\dim}$. 
Note that $\widehat{\P}\n$ does not depend on the current, meaning the current is not coupled across interior faces in its diagonal block, allowing the efficient solution procedure outlined above. 
Here, $\kappa\geq 0$ is a free parameter that we elect to set to $\alpha/2$ to exactly match the analogous term in the \Pone LO system \eqref{eq:p1-m0}. 
The numerical fluxes on the boundary are the same as for \Pone from \eqref{eq:P1_Jhat_bdr} and \eqref{eq:P1_Phat_bdr}. 

Observe that the LDG numerical fluxes add $s/2\jump{\vec{J}\cdot\n}$ to the \Pone zeroth moment numerical flux in \eqref{eq:P1_Jhat} and that $1/6\alpha\jump{\vec{J}\cdot\n}$ is replaced with $-s/2\jump{\varphi}$ in the first moment numerical flux from \eqref{eq:P1_Phat}. 
The LDG SMM numerical fluxes are then 
	\begin{subequations}
	\begin{equation}
		\widehat{\vec{J}}_\text{LDG}\cdot\n = \avg{\vec{J}\cdot\n} + \frac{s}{2}\jump{\vec{J}\cdot\n} + \frac{\alpha}{2}\jump{\varphi} + \frac{1}{2}\jump{\beta - s\Jho\cdot\n} \,, 
	\end{equation}
	\begin{equation}
		\widehat{\P}_\text{LDG}\n = \frac{\n}{3}\!\paren{\avg{\varphi} - \frac{s}{2}\jump{\varphi}} + \avg{\T\n} + \frac{1}{2}\jump{\Phoplus - \Phominus + \frac{\n}{3}s\phiho} \,,
	\end{equation}
	\end{subequations}
for interior faces and 
	\begin{subequations}
	\begin{equation}
		\widehat{\vec{J}}_\text{LDG}\cdot\n = \frac{\vec{J}\cdot\n}{2} + \frac{\alpha}{2}\varphi + \frac{1}{2}\beta + J_\text{in} \,, 
	\end{equation}
	\begin{equation}
		\widehat{\P}_\mathrm{LDG}\n = \Phoplus + \vec{P}_\text{in} + \frac{\n}{6\alpha}\vec{J}\cdot\n+ \frac{\n}{6}\varphi - \frac{\n}{6\alpha}\Jho\cdot\n - \frac{\n}{6}\phiho\,,
	\end{equation}
	\end{subequations}
for boundary faces. 
As with \Pone, the pressure term is closed with \eqref{eq:SMM_P}. 
Replacing the numerical fluxes in the \Pone weak form \eqref{eq:p1-weak}, the consistent LDG LO problem is: find $(\varphi,\vec{J}) \in Y_1 \times [Y_1]^{\dim}$ such that 
	\begin{subequations}
	\begin{multline}
		\int_{\Gamma_0} \jump{u}\paren{\avg{\vec{J}\cdot\n} + \frac{s}{2}\jump{\vec{J}\cdot\n}} \ud s + \frac{\alpha}{2}\int_{\Gamma_0} \jump{u}\jump{\varphi} \ud s + \frac{1}{2}\int_{\Gamma_b} u\, \vec{J}\cdot\n \ud s + \frac{\alpha}{2}\int_{\Gamma_b} u\varphi \ud s \\
		- \int \nablah u \cdot\vec{J} \ud \x + \int \sigma_a\, u \varphi \ud \x = \int u\, Q_0 \ud \x - \int_{\Gamma_b} u\, J_\text{in} \ud s + \mathcal{R}_\text{LDG}^0(u) \,, \quad \forall u \in Y_1 \,,
	\end{multline}
	\begin{multline}
		\frac{1}{3}\int_{\Gamma_0} \jump{\vec{v}\cdot\n}\paren{\avg{\varphi} - \frac{s}{2}\jump{\varphi}} \ud s + \frac{1}{6\alpha}\int_{\Gamma_b} (\vec{v}\cdot\n)(\vec{J}\cdot\n) \ud s + \frac{1}{6}\int_{\Gamma_b} \vec{v}\cdot\n\, \varphi \ud s \\ - \frac{1}{3}\int \nablah\cdot\vec{v}\, \varphi \ud \x + \int \sigma_t\, \vec{v}\cdot\vec{J} \ud \x = \int \vec{v}\cdot\vec{Q}_1 \ud \x - \int \vec{v}\cdot\vec{P}_\text{in} \ud s + \mathcal{R}_\text{LDG}^1(\vec{v}) \,, \quad \forall \vec{v} \in [Y_1]^{\dim} \,. 
	\end{multline}
	\end{subequations}
The LDG correction source terms are 
	\begin{subequations}
	\begin{equation} \label{eq:LDG_R0}
		\mathcal{R}_\text{LDG}^0(u) = -\frac{1}{2}\int_{\Gamma_0} \jump{u} \jump{\beta} - \frac{1}{2}\int_{\Gamma_b} u\, \beta \ud s + \frac{s}{2}\int_{\Gamma_0} \jump{u} \jump{\Jho\cdot\n} \ud s \,,
	\end{equation}
	\begin{multline} \label{eq:LDG_R1}
		\mathcal{R}_\text{LDG}^1(\vec{v}) = -\int_{\Gamma_0} \jump{\vec{v}}\cdot\avg{\T\n} \ud s - \frac{1}{2}\int_{\Gamma_0} \jump{\vec{v}}\cdot\jump{\Phoplus - \Phominus + \frac{\n}{3}s\phiho} \ud s \\- \int_{\Gamma_b} \vec{v} \cdot \paren{\Phoplus - \frac{\n}{6\alpha}\Jho\cdot\n - \frac{\n}{6}\phiho} \ud s + \int \nablah\vec{v} : \T \ud \x \,. 
	\end{multline}
	\end{subequations}

\subsubsection{Interior Penalty}
The IP numerical fluxes are 
	\begin{subequations}
	\begin{equation} \label{eq:ip_J}
		\widehat{\vec{J}}\cdot\n = \avg{\vec{J}\cdot\n} + \kappa \jump{\varphi} \,, 
	\end{equation}
	\begin{equation} \label{eq:ip_P}
		\widehat{\P}\n = \frac{\n}{3}\avg{\varphi} \,, 
	\end{equation}
	\end{subequations}
where the parameter $\kappa \propto \sigma_t^{-1} h^{-1}$ is required to ensure stability as $h\rightarrow 0$ and the proportionality constant is typically problem dependent and requires tuning \cite{Arnold2002}. 
Compared to the zeroth moment \Pone numerical flux in \eqref{eq:P1_Jhat}, \eqref{eq:ip_J} differs in the choice of coefficient in the jump-jump bilinear form, $\kappa$ versus $\alpha/2$. 
For the first moment, \eqref{eq:ip_P} has removed the term depending on $\jump{\vec{J}\cdot\n}$ in \eqref{eq:P1_P}. 
The interior numerical fluxes are then 
	\begin{subequations}
	\begin{equation}
		\widehat{\vec{J}}_\text{IP}\cdot\n = \avg{\vec{J}\cdot\n} + \kappa \jump{\varphi} + \frac{1}{2}\jump{\beta} + \paren{\frac{\alpha}{2} - \kappa} \jump{\phiho} \,,
	\end{equation}
	\begin{equation}
		\widehat{\P}_\text{IP}\n = \frac{\n}{3}\avg{\varphi} + \avg{\T\n} + \frac{1}{2}\jump{\Phoplus - \Phominus} \,. 
	\end{equation}
	\end{subequations}
The factor of $\alpha/2 - \kappa$ occurs because $\beta/2$ implicitly includes $-\alpha/2 \phiho$. 
We present results with $\kappa$ prescribed according to the modified interior penalty formulation from \cite{WR} defined as
	\begin{equation}
		\kappa_\text{MIP} = \max\paren{\kappa_\text{IP}\,, \alpha/2} \,,
	\end{equation}
where $\kappa_\text{IP} = \frac{C}{3\sigma_t h}$ is the standard interior penalty parameter with problem-dependent constant, $C$. 
In the under-resolved thick diffusion limit, $\kappa_\text{MIP} \rightarrow \alpha/2$, causing $(\alpha/2-\kappa)\jump{\phiho}$ to vanish. 
The IP LO system is: find $(\varphi,\vec{J}) \in Y_1 \times [Y_1]^{\dim}$ such that 
	\begin{subequations}
	\begin{multline}
		\int_{\Gamma_0} \jump{u}\avg{\vec{J}\cdot\n} \ud s + \int_{\Gamma_0} \kappa\,\jump{u}\jump{\varphi} \ud s + \frac{1}{2}\int_{\Gamma_b} u\, \vec{J}\cdot\n \ud s + \frac{\alpha}{2}\int_{\Gamma_b} u\varphi \ud s \\
		- \int \nablah u \cdot\vec{J} \ud \x + \int \sigma_a\, u \varphi \ud \x = \int u\, Q_0 \ud \x - \int_{\Gamma_b} u\, J_\text{in} \ud s + \mathcal{R}_\text{IP}^0(u) \,, \quad \forall u \in Y_1 \,,
	\end{multline}
	\begin{multline}
		\frac{1}{3}\int_{\Gamma_0} \jump{\vec{v}\cdot\n}\avg{\varphi} \ud s + \frac{1}{6\alpha}\int_{\Gamma_b} (\vec{v}\cdot\n)(\vec{J}\cdot\n) \ud s + \frac{1}{6}\int_{\Gamma_b} \vec{v}\cdot\n\, \varphi \ud s \\ - \frac{1}{3}\int \nablah\cdot\vec{v}\, \varphi \ud \x + \int \sigma_t\, \vec{v}\cdot\vec{J} \ud \x = \int \vec{v}\cdot\vec{Q}_1 \ud \x - \int \vec{v}\cdot\vec{P}_\text{in} \ud s + \mathcal{R}_\text{IP}^1(\vec{v}) \,, \quad \forall \vec{v} \in [Y_1]^{\dim} \,, 
	\end{multline}
	\end{subequations}
where 
	\begin{subequations}
	\begin{equation}
		\mathcal{R}_\text{IP}^0(u) = -\frac{1}{2}\int_{\Gamma_0} \jump{u} \jump{\beta} - \frac{1}{2}\int_{\Gamma_b} u\, \beta \ud s + \int_{\Gamma_0} \paren{\kappa - \frac{\alpha}{2}}\jump{u} \jump{\phiho} \ud s \,,
	\end{equation}
	\begin{multline}
		\mathcal{R}_\text{IP}^1(\vec{v}) = -\int_{\Gamma_0} \jump{\vec{v}}\cdot\avg{\T\n} \ud s - \frac{1}{2}\int_{\Gamma_0} \jump{\vec{v}}\cdot\jump{\Phoplus - \Phominus} \ud s \\- \int_{\Gamma_b} \vec{v} \cdot \paren{\Phoplus - \frac{\n}{6\alpha}\Jho\cdot\n - \frac{\n}{6}\phiho} \ud s + \int \nablah\vec{v} : \T \ud \x \,.
	\end{multline}
	\end{subequations}

\subsubsection{Alternative Boundary Conditions} \label{sec:alt_bcs}
The boundary conditions derived from the moments of DG \Sn transport have a different structure than those used for the independent DG SMMs derived in \cite{olivier_smm}. 
We refer to these two types of boundary conditions as half and full-range, respectively, owing to the domain of integration associated with the boundary corrections. 
It was found that the boundary conditions affect the solvability of the LO system and the performance of the consistent SMMs in the thick diffusion limit. 
The boundary conditions from \cite{olivier_smm} are of the form
	\begin{subequations}
	\begin{equation}
		\widehat{\vec{J}}\cdot\n = \alpha \varphi + 2 J_\text{in} \,,
	\end{equation}
	\begin{equation}
		\widehat{\P}\n = \frac{\n}{3}\varphi \,. 
	\end{equation}
	\end{subequations}
Comparing these boundary conditions to the DG \Sn moment boundary conditions from \eqref{eq:dgsn_Jflux_bdr} and \eqref{eq:dgsn_Pflux_bdr}, alternative boundary fluxes for LDG and IP are
	\begin{subequations}
	\begin{equation}
	\begin{aligned}
		\widehat{\vec{J}}_\text{full}\cdot\n &= \widehat{\vec{J}}_{\smallHO}\cdot\n + \alpha \varphi + 2J_\text{in} - \paren{\alpha \phiho + 2J_\text{in}} \\
		&= \alpha \varphi + \Jhoplus - \alpha \phiho + J_\text{in} \,, 
	\end{aligned}
	\end{equation}
	\begin{equation}
	\begin{aligned}
		\widehat{\P}_\text{full}\n &= \widehat{\P}_{\smallHO}\n + \frac{\n}{3}\varphi - \frac{\n}{3}\phiho \\
		&= \frac{\n}{3}\varphi + \Phoplus - \frac{\n}{3}\phiho + \vec{P}_\text{in} \,. 
	\end{aligned}
	\end{equation}
	\end{subequations}
With these full range boundary fluxes, the LDG LO problem is: find $(\varphi,\vec{J}) \in Y_1\times [Y_1]^{\dim}$ such that 
\begin{subequations}
	\begin{multline}
		\int_{\Gamma_0} \jump{u}\paren{\avg{\vec{J}\cdot\n} + \frac{s}{2}\jump{\vec{J}\cdot\n}} \ud s + \frac{\alpha}{2}\int_{\Gamma_0} \jump{u}\jump{\varphi} \ud s + \alpha\int_{\Gamma_b} u\varphi \ud s \\
		- \int \nablah u \cdot\vec{J} \ud \x + \int \sigma_a\, u \varphi \ud \x = \int u\, Q_0 \ud \x - 2\int_{\Gamma_b} u\, J_\text{in} \ud s + \mathcal{R}_\mathrm{LDG,full}^0(u) \,, \quad \forall u \in Y_1 \,,
	\end{multline}
	\begin{multline}
		\frac{1}{3}\int_{\Gamma_0} \jump{\vec{v}\cdot\n}\paren{\avg{\varphi} - \frac{s}{2}\jump{\varphi}} \ud s + \frac{1}{3}\int_{\Gamma_b} \vec{v}\cdot\n\, \varphi \ud s \\ - \frac{1}{3}\int \nablah\cdot\vec{v}\, \varphi \ud \x + \int \sigma_t\, \vec{v}\cdot\vec{J} \ud \x = \int \vec{v}\cdot\vec{Q}_1 \ud \x - \int \vec{v}\cdot\vec{P}_\text{in} \ud s + \mathcal{R}_\mathrm{LDG,full}^1(\vec{v}) \,, \quad \forall \vec{v} \in [Y_1]^{\dim} \,, 
	\end{multline}
	\end{subequations}
where 
	\begin{subequations}
	\begin{equation}
		\mathcal{R}_\mathrm{LDG,full}^0(u) = -\frac{1}{2}\int_{\Gamma_0} \jump{u}\jump{\beta} \ud s + \frac{s}{2} \int_{\Gamma_0} \jump{u}\jump{\Jho\cdot\n} \ud s - \int_{\Gamma_b} u\paren{\Jhoplus - \alpha \phiho - J_\text{in}} \ud s \,,
	\end{equation}
	\begin{multline}
		\mathcal{R}_\mathrm{LDG,full}^1(\vec{v}) = -\int_{\Gamma_0} \jump{\vec{v}}\cdot\avg{\T\n} \ud s - \frac{1}{2}\int_{\Gamma_0} \jump{\vec{v}}\cdot\jump{\Phoplus - \Phominus + \frac{\n}{3}s\phiho} \ud s \\- \int_{\Gamma_b} \vec{v} \cdot \paren{\Phoplus - \frac{\n}{3}\phiho} \ud s + \int \nablah\vec{v} : \T \ud \x \,. 
	\end{multline}
	\end{subequations}
With full range boundary conditions, the IP LO problem is: find $(\varphi,\vec{J}) \in Y_1 \times [Y_1]^{\dim}$ such that 
	\begin{subequations} \label{eq:iplo_full}
	\begin{multline}
		\int_{\Gamma_0} \jump{u}\avg{\vec{J}\cdot\n} \ud s + \int_{\Gamma_0} \kappa\,\jump{u}\jump{\varphi} \ud s + \alpha\int_{\Gamma_b} u\varphi \ud s \\
		- \int \nablah u \cdot\vec{J} \ud \x + \int \sigma_a\, u \varphi \ud \x = \int u\, Q_0 \ud \x - 2\int_{\Gamma_b} u\, J_\text{in} \ud s + \mathcal{R}_\mathrm{IP,full}^0(u) \,, \quad \forall u \in Y_1 \,,
	\end{multline}
	\begin{multline}
		\frac{1}{3}\int_{\Gamma_0} \jump{\vec{v}\cdot\n}\avg{\varphi} \ud s + \frac{1}{3}\int_{\Gamma_b} \vec{v}\cdot\n\, \varphi \ud s - \frac{1}{3}\int \nablah\cdot\vec{v}\, \varphi \ud \x \\+ \int \sigma_t\, \vec{v}\cdot\vec{J} \ud \x = \int \vec{v}\cdot\vec{Q}_1 \ud \x - \int \vec{v}\cdot\vec{P}_\text{in} \ud s + \mathcal{R}_\mathrm{IP,full}^1(\vec{v}) \,, \quad \forall \vec{v} \in [Y_1]^{\dim} \,,
	\end{multline}
	\end{subequations}
where 
	\begin{subequations}
	\begin{equation}
		\mathcal{R}_\mathrm{IP,full}^0(u) = -\frac{1}{2}\int_{\Gamma_0} \jump{u} \jump{\beta} \ud s - \int_{\Gamma_b} u\paren{\Jhoplus - \alpha \phiho - J_\text{in}} \ud s + \int_{\Gamma_0} \paren{\kappa - \frac{\alpha}{2}}\, \jump{u} \jump{\phiho} \ud s \,,
	\end{equation}
	\begin{multline}
		\mathcal{R}_\mathrm{IP,full}^1(\vec{v}) = -\int_{\Gamma_0} \jump{\vec{v}}\cdot\avg{\T\n} \ud s - \frac{1}{2}\int_{\Gamma_0} \jump{\vec{v}}\cdot\jump{\Phoplus - \Phominus} \ud s \\- \int_{\Gamma_b} \vec{v} \cdot \paren{\Phoplus - \frac{\n}{3}\phiho} \ud s + \int \nablah\vec{v} : \T \ud \x \,. 
	\end{multline}
	\end{subequations}

\subsection{Analysis of Consistent Closures}
The correction source terms discussed above contain model and discretization correction terms. 
We define discretization corrections as terms that are on the order of the discretization error and can be removed without degrading the physics fidelity of the algorithm resulting in an independent SMM. 
Since the approximation spaces for the LO solution and the moments of the HO solution are the same, the volumetric terms from the LO system exactly match the volumetric terms in the moments of the HO system. 
Thus, discretization corrections arise only on the interior and boundary faces of the mesh. 
Note that, when the true solution is continuous, $\jump{\psi_d}$ is on the order of the discretization error, hence, terms in the correction sources that depend on a jump in a HO variable can be ignored without altering the overall algorithm's spatial convergence rate. 

Consider, for example, the LDG correction source terms \eqref{eq:LDG_R0} and \eqref{eq:LDG_R1}. 
The interior face terms depending on $\jump{\beta}$, $\jump{\Jho\cdot\n}$, $\jump{\phiho}$, and $\jump{\Phoplus - \Phominus}$ can all be dropped and, by our definition, are discretization corrections. 
In contrast, the terms 
	\begin{equation}
		\int_{\Gamma_0} \jump{\vec{v}}\cdot\avg{\T\n} \ud s - \int \nablah\vec{v} : \T \ud \x 
	\end{equation}
vanish only when the solution is linearly anisotropic, which we consider to be model correction terms. 
In practice, we have seen that neglecting the discretization correction terms results in lower solution quality, especially for the current. 
Because the correction sources in a consistent method contain discretization corrections in addition to the model corrections required by an independent SMM, we postulate that consistent methods will require more fixed-point iterations than independent methods. 
However, these additional terms are on the order of the discretization error so that consistent and independent methods should converge at similar rates on a mesh that is refined enough. 
Finally, note that the independent IP SMM in \cite{olivier_smm} is equivalent to the consistent IP SMM with full-range boundary conditions from \eqref{eq:iplo_full} when the discretization corrections are neglected. 

\section{Numerical Results}\label{sec:results}
We now present numerical results demonstrating the accuracy and performance of the methods. 
The implementations used the MFEM finite element framework \cite{mfem-paper} and and were solved with the \texttt{BoomerAMG} preconditioner from the \textit{hypre} sparse linear algebra library \cite{hypre} and the \texttt{KINSOL} fixed-point solver with Anderson acceleration from the Sundials package \cite{hindmarsh2005sundials}. 
The \Pone discretization is solved with the distributed version of SuperLU \cite{lidemmel03}. 

\subsection{Accuracy}
\begin{figure}
\centering
\begin{subfigure}{.49\textwidth}
	\centering
	\includegraphics[width=\textwidth]{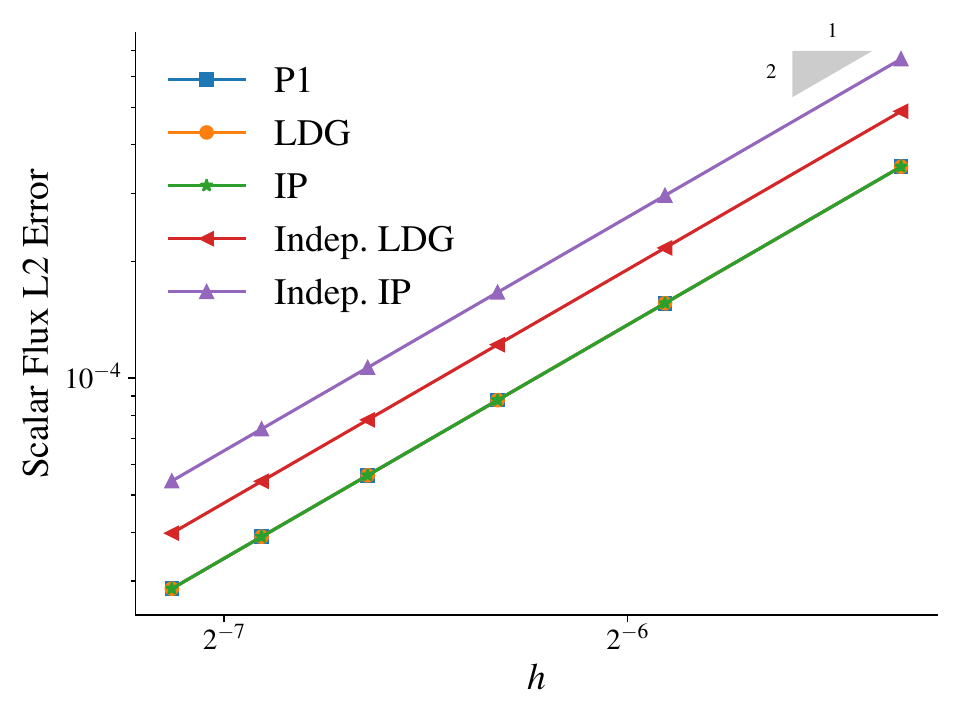}
	\caption{Scalar Flux}
	\label{fig:mms_phi}
\end{subfigure}
\begin{subfigure}{.49\textwidth}
	\centering
	\includegraphics[width=\textwidth]{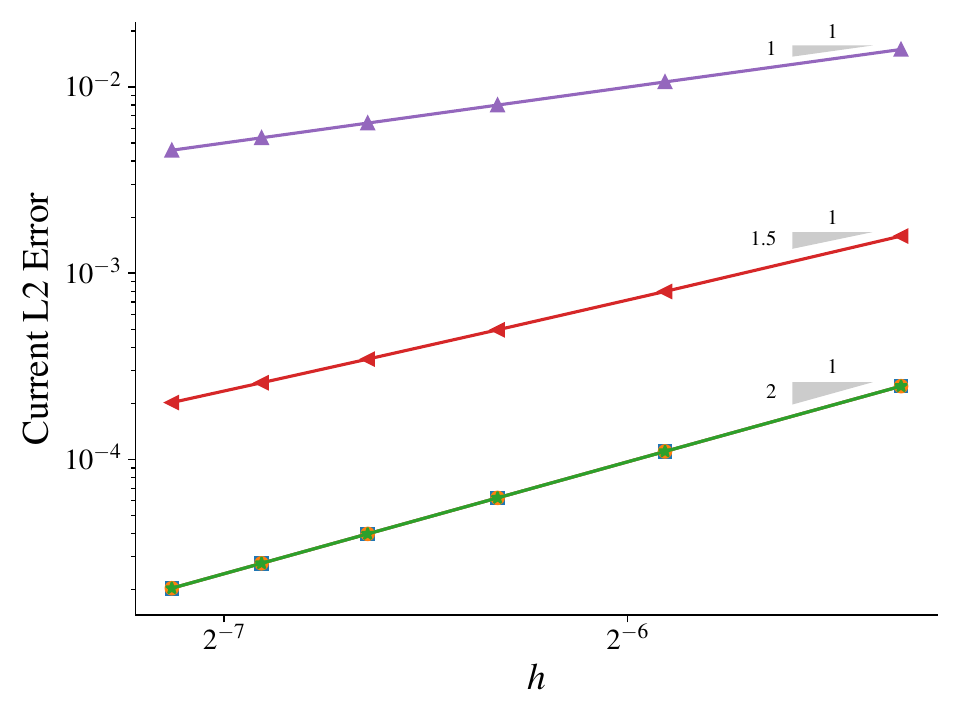}
	\caption{Current}
	\label{fig:mms_J}
\end{subfigure}
\caption{Errors on the MMS problem. The \Pone, LDG, and IP methods were equivalent to the iterative tolerance of $10^{-10}$.}
\end{figure}

The accuracy of the schemes are investigated with the Method of Manufactured Solutions (MMS).
The computational domain is $\D = [0,1]^2$. 
The solution is set to 
	\begin{multline} \label{eq:psi_mms}
		\psi = \frac{1}{4\pi}\bigr(\sin(\pi x) \sin(\pi y) + (\Omegahat_x + \Omegahat_y)\sin(2\pi x) \sin(2\pi y)/2 \\+ (\Omegahat_x^2 + \Omegahat_y^2)\sin\!\paren{\frac{3\pi(x+\delta)}{1+2\delta}}\sin\!\paren{\frac{3\pi(y+\delta)}{1+2\delta}}/4 + 2\bigr) \,, 
	\end{multline}
where $\delta = 0.05$ allows testing the accuracy of the methods on problems with spatially dependent, anisotropic inflow boundary conditions. 
The scalar flux and current solutions are then: 
	\begin{subequations}
	\begin{equation} \label{eq:phi_mms}
		\varphi = \sin(\pi x) \sin(\pi y) + \frac{1}{6}\sin\!\paren{\frac{3\pi(x+\delta)}{1+2\delta}}\sin\!\paren{\frac{3\pi(y+\delta)}{1+2\delta}} + 2 \,, 
	\end{equation}
	\begin{equation} \label{eq:J_mms}
		\vec{J} = \frac{1}{6}\sin(2\pi x) \sin(2\pi y) \begin{bmatrix} 1 \\ 1 \end{bmatrix} \,. 
	\end{equation}
	\end{subequations}
The material data are $\sigma_t = \SI{2}{\per\cm}$ and $\sigma_s = \SI{1.9}{\per\cm}$. 
The fixed-source MMS problem is solved with fixed-point iteration to a tolerance of $10^{-10}$ with level symmetric $S_4$ quadrature. 
The error is computed in the $L^2(\D)$ and $[L^2(\D)]^{\dim}$ for the scalar flux and current, respectively. 
Figures \ref{fig:mms_phi} and \ref{fig:mms_J} show the errors as the mesh is refined, comparing the consistent SMMs and the independent variants of LDG and IP. 
For the scalar flux, all methods produce the optimal second-order convergence. 
However, for the current, only the consistent methods are able to achieve optimal accuracy with the independent LDG method converging the current at $\mathcal{O}(h^{3/2})$ and independent IP at first-order. 
Notably, all the consistent methods produce the same error for both the scalar flux and current. 

\begin{figure}
\centering
\begin{subfigure}{.49\textwidth}
	\centering
	\includegraphics[width=\textwidth]{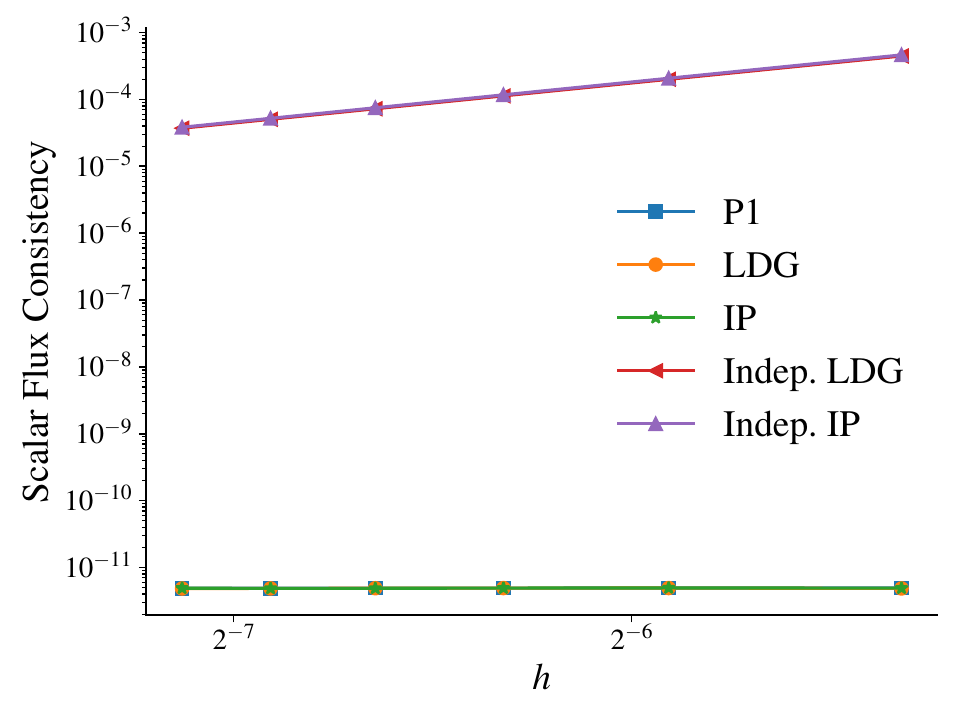}
	\caption{Scalar Flux}
	\label{fig:consistency_phi}
\end{subfigure}
\begin{subfigure}{.49\textwidth}
	\centering
	\includegraphics[width=\textwidth]{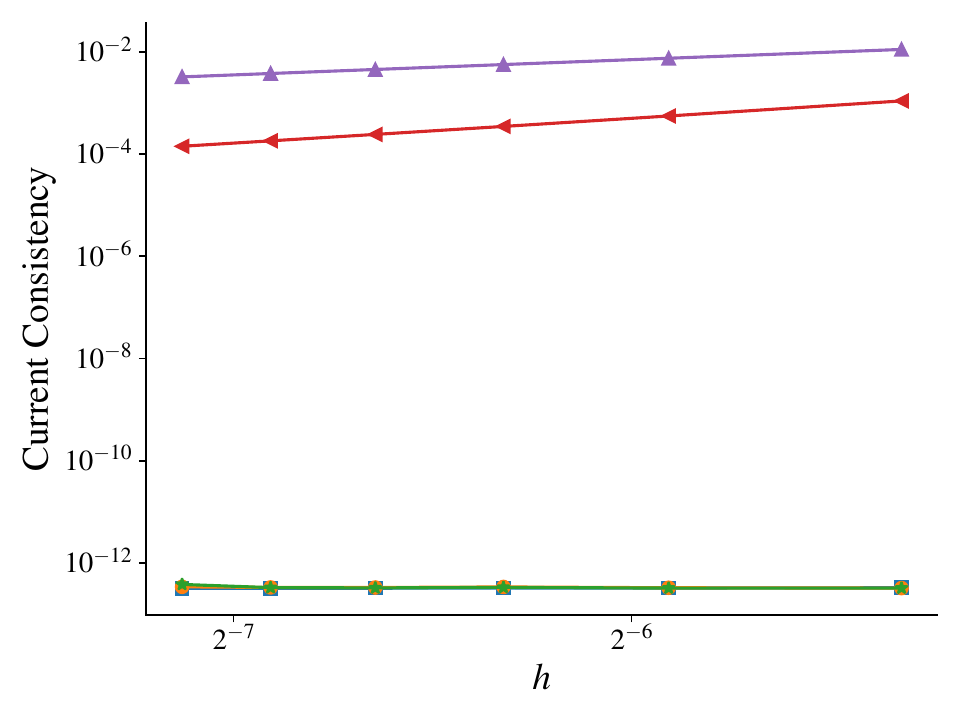}
	\caption{Current}
	\label{fig:consistency_J}
\end{subfigure}
\caption{Consistency between the HO and LO solution variables on the MMS problem. The \Pone, LDG, and IP methods were equivalent to below the iterative tolerance of $10^{-10}$.}
\end{figure}

Next, we demonstrate that the consistent methods produce HO and LO solutions that match up to iterative tolerances independent of the mesh size. 
Figures \ref{fig:consistency_phi} and \ref{fig:consistency_J} show the difference between the LO solution and the zeroth and first angular moments of the HO solution. 
In other words, these figures plot 
	\begin{equation}
		\| \varphi - \phiho \| \,, \quad \| \vec{J} - \Jho \| \,, 
	\end{equation}
where $\phi_\text{HO} = \sum_d w_d\,\psi_d$ and $\vec{J}_\text{HO} = \sum_d w_d\,\Omegahat_d\,\psi_d$. 
Observe that the three consistent schemes produce consistency errors below the iterative tolerance of $10^{-10}$ and that the error is independent of the mesh size. 
By contrast, the independent methods have consistency error on the order of the LO discretization error and thus the consistency error depends on the mesh size. 

\subsection{Thick Diffusion Limit}
We now compare the iterative efficiency of the methods in the asymptotic, thick diffusion limit. 
The domain is $[0,1]^2$. 
The materials are set to 
	\begin{equation}
		\sigma_t = 1/\epsilon \,, \quad \sigma_s = 1/\epsilon - \epsilon \,, \quad q = \epsilon \,, 
	\end{equation}
where $\epsilon \in (0,1]$ such that $\epsilon \rightarrow 0$ induces the diffusion limit. 
Level symmetric $S_4$ angular quadrature is used with a coarse $8\times 8$ mesh. 
This is a numerically challenging, but common in practice, regime where robust performance is crucial. 

\begin{figure}
\centering
\begin{subfigure}{.32\textwidth}
	\centering
	\includegraphics[width=\textwidth]{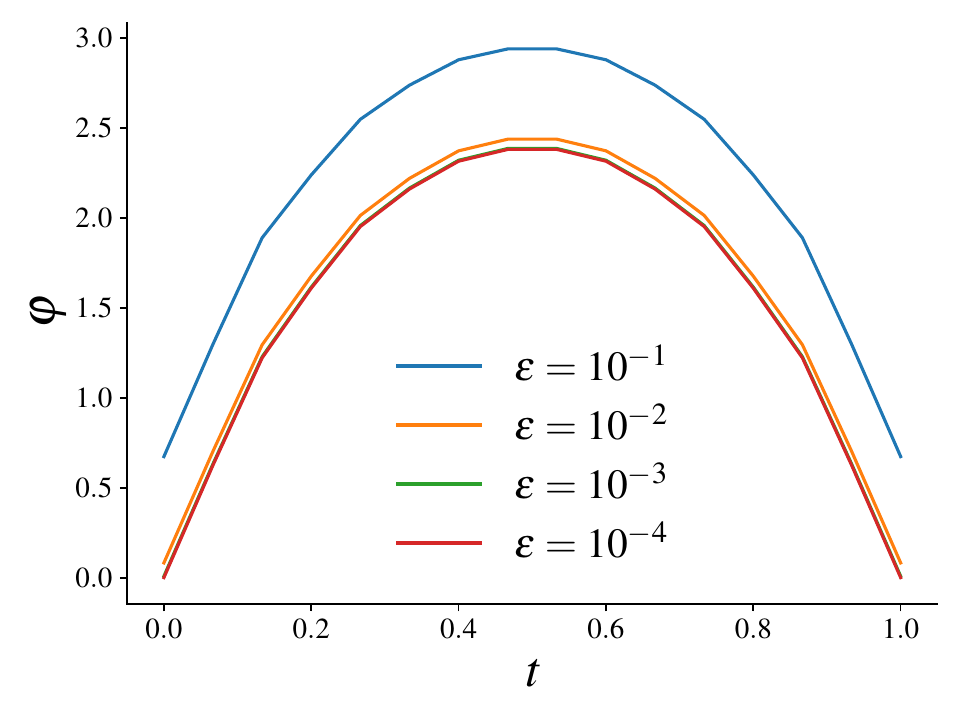}
	\caption{\Pone}
\end{subfigure}
\begin{subfigure}{.32\textwidth}
	\centering
	\includegraphics[width=\textwidth]{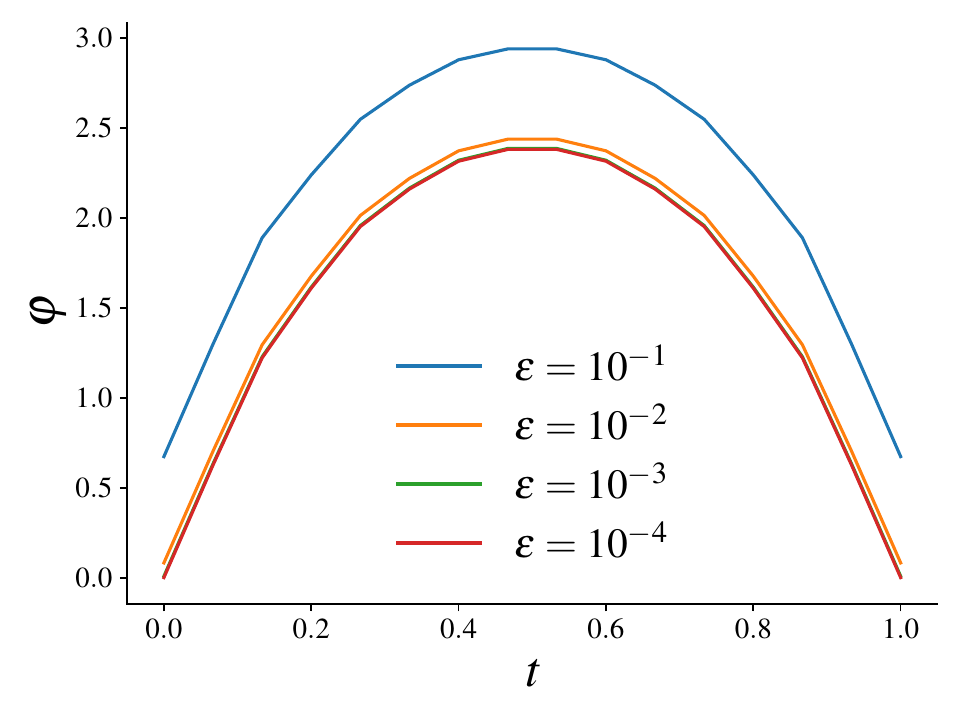}
	\caption{LDG}
\end{subfigure}
\begin{subfigure}{.32\textwidth}
	\centering
	\includegraphics[width=\textwidth]{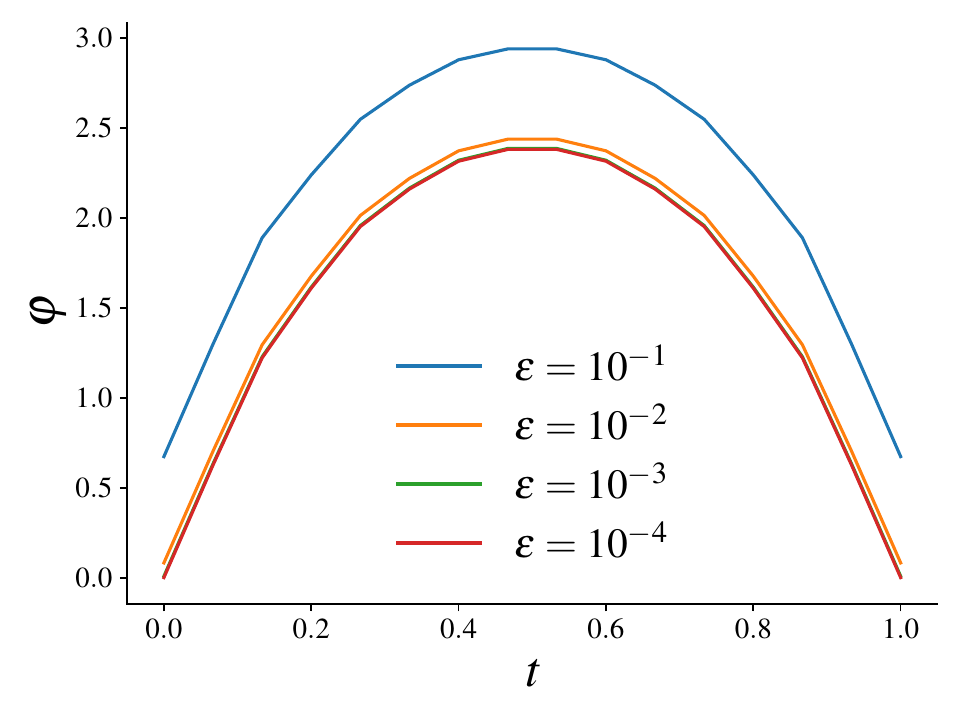}
	\caption{IP}
\end{subfigure}
\caption{Lineouts of the scalar flux along $y=\SI{0.5}{\cm}$ for the consistent SMMs as $\epsilon \rightarrow 0$. }
\label{fig:eps_line}
\end{figure}

\begin{table}
\centering
\caption{Iterations for SMM to converge in thick diffusion limit.}
\label{tab:tdl}
\begin{tabular}{ccccccccccc}
\toprule
 &  & \multicolumn{3}{c}{LDG}  &  & \multicolumn{4}{c}{IP} \\
\cmidrule{3-5}\cmidrule{7-10}
$\epsilon$ & P1 & Full & Half & Indep. & & Full & Half & Unmodified & Indep. \\
\midrule
$10^{-1}$ & 9 & 13 & 9 & 9 & & 13 & 9 & 13 & 9 \\
$10^{-2}$ & 5 & 16 & 5 & 5 & & 16 & 5 & 54 & 5 \\
$10^{-3}$ & 4 & 15 & 3 & 3 & & 15 & 3 & 151 & 3 \\
$10^{-4}$ & 3 & 11 & 3 & 3 & & 11 & 3 & 70 & 3 \\
\bottomrule
\end{tabular}
\end{table}
Table \ref{tab:tdl} shows the number of fixed-point iterations for convergence to a tolerance of $10^{-6}$ as $\epsilon\rightarrow 0$. 
The consistent $P_1$, LDG, and IP SMMs are compared to the independent variants of LDG and IP. 
The LDG and IP methods are compared with both full and half-range boundary conditions. 
In addition, the consistent IP scheme without the modified interior penalty parameter specification from \cite{WR} is shown. 
Aside from the unmodified IP method, all schemes achieve iterative efficiency independent of $\epsilon$. 
Convergence of the consistent LDG and IP methods with half-range boundary conditions is comparable to that of the P$_1$ method and convergence is significantly improved compared to the full-range boundary conditions.
The IP method without the modified penalty parameter scales poorly with $\epsilon$. 
This problem is severely under-resolved, in which case the IP discretization corrections are relatively large compared to the physics corrections, slowing iterative convergence.
The modified penalty parameter ensures that the discretization corrections remain bounded as $\epsilon\rightarrow 0$. 
Lineouts of the solutions for the consistent P$_1$, LDG, and IP SMMs are shown in Fig.~\ref{fig:eps_line}, illustrating that physically realistic solutions are obtained as $\epsilon \rightarrow 0$.

\subsection{Crooked Pipe Problem}
\begin{figure}
\centering
\includegraphics[width=.65\textwidth]{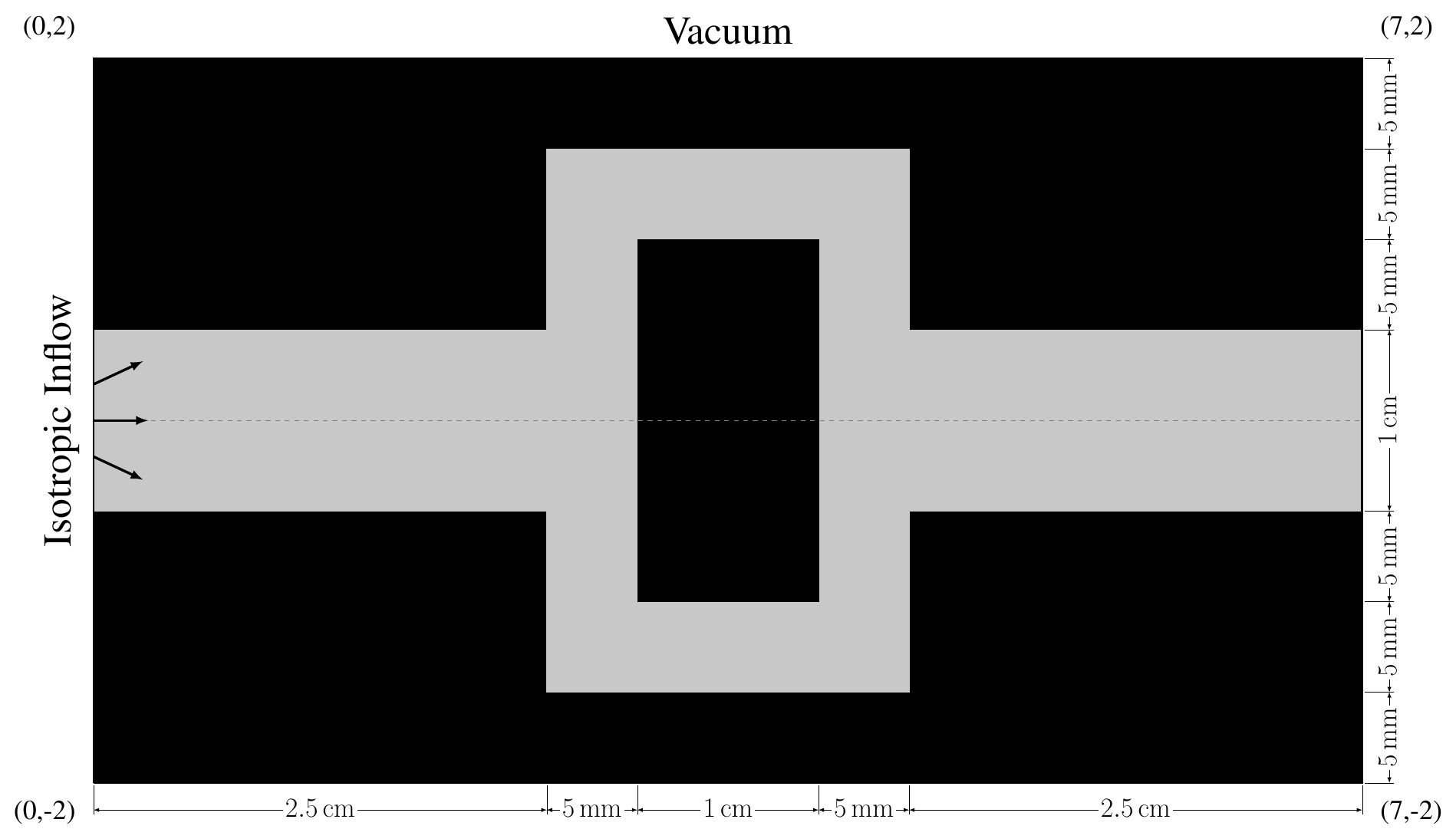}
\caption{Depiction of the geometry of the crooked pipe benchmark problem. The optically thin pipe and optically thick wall are depicted as gray and black, respectively.}
\label{fig:cp_geom}
\end{figure}
Here, we compare performance in parallel and solution quality of the consistent and independent SMMs on a steady-state, linear transport variant of the crooked pipe benchmark. 
The geometry, depicted in Fig.~\ref{fig:cp_geom}, consists of two materials, the optically thick wall and the optically thin pipe, which have a $1000\times$ difference in total cross section between the two regions. 
The wall and pipe are characterized by 
	\begin{subequations}
	\begin{align}
		\textnormal{Wall:}\hspace{4ex}\sigma_t & = \SI{200}{\per\cm}\,, \quad \sigma_a = 10^{-3}\si{\per\cm} \,, \quad q = 10^{-7}\si{\per\cm\cubed\per\str\per\s}\,, \\
		\textnormal{Pipe:}\hspace{4ex}\sigma_t & = \SI{0.2}{\per\cm}\,, \quad \sigma_a = 10^{-3}\si{\per\cm} \,, \quad q = 10^{-7}\si{\per\cm\cubed\per\str\per\s}\,. 
	\end{align}
	\end{subequations}
The absorption and source are chosen to correspond to the first backward Euler time step of a thermal radiative transfer calculation with $c\Delta t = 10^3\si{\cm}$ and initial condition $\psi_0 = 10^{-4} \si{\per\cm\squared\per\str\per\s}$. 
The isotropic, inflow boundary condition is 
	\begin{equation}
		\bar{\psi} = \begin{cases}
			\frac{1}{2\pi} \,, & x=0 \ \mathrm{and} \ |y| \leq \SI{0.5}{\cm} \\ 
			0 \,, & \mathrm{otherwise}
		\end{cases} \,, 
	\end{equation}
so that radiation isotropically enters the left entrance of the pipe with vacuum boundary conditions applied elsewhere. 
A reflection plane along the line $y=0$ is used to halve the computational domain.
We use level symmetric $S_{12}$ angular quadrature. 

\subsubsection{Parallel Performance}
The efficiency of the outer fixed-point iteration with preconditioned, linear inner iteration is investigated under weak scaling wherein the MPI-decomposed mesh is uniformly refined as the number of processors is increased so that the work per processor remains constant. 
The base mesh has $224 \times 64$ uniform quadrilateral elements that are aligned with the materials in the problem. 
Across four refinements of this base mesh, the mesh size ranges from $h = \SI{3.125e-2}{\cm}$ to $h = \SI{1.953e-3}{\cm}$. 
We compare fixed-point iteration and Anderson-accelerated fixed-point iteration \cite{10.1145/321296.321305} with five Anderson vectors.
The outer tolerance is $10^{-6}$. 
The P$_1$ LO system is solved with SuperLU. 
All other methods are solved with conjugate gradient preconditioned with one V-cycle of AMG to a tolerance of $10^{-8}$. 
The solution from the previous outer iteration is used as an initial guess for the inner iteration so that the initial guess becomes progressively more accurate as the outer iteration converges. 
The streaming and collision operator is inverted directly with a full parallel transport sweep. 
The sweep is ordered so that the reflection boundary condition is inverted exactly. 
The calculations were run on 22-core Intel Skylake Gold processors with two sockets and 376 GB of memory per node. 
The runtime is reported as the minimum of three repeated executions for each method and refinement. 

The number of outer fixed-point iterations is shown in Table \ref{tab:fp_outer} for the consistent P$_1$, LDG, and IP and independent LDG and IP methods.
The final \Pone data point could not be obtained due to the SuperLU solve of the \Pone LO system exceeding the memory available on six nodes. 
For P$_1$ and the independent methods, the number of outer iterations increases as the mesh is refined. 
In contrast, the consistent LDG and IP methods show decreasing iterations as the mesh is refined. 
This could be explained by the discretization corrections becoming smaller in magnitude as the mesh is refined and the discretization error is reduced. 
Unlike on the single material thick diffusion limit problem, full or half-range boundary conditions did not significantly affect convergence. 
For the least resolved mesh, the independent SMMs converged 1.7x faster than P$_1$ and 2.5x faster than the LDG and IP consistent methods. 
However, as the mesh is refined all the methods eventually converge at the same rate. 
\begin{table}[h!]
\centering
\caption{Number of outer fixed-point iterations for the crooked pipe problem.}
\begin{tabular}{ccccccccccc}
\toprule
 &  &  & \multicolumn{3}{c}{LDG}  &  & \multicolumn{3}{c}{IP} \\
\cmidrule{4-6}\cmidrule{8-10}
$N_e$ & Proc. & P1 & Full & Half & Indep. & & Full & Half & Indep. \\
\midrule
\num{14336} & 1 & 68 & 96 & 97 & 38 & & 109 & 108 & 38 \\
\num{57344} & 4 & 73 & 90 & 89 & 47 & & 98 & 98 & 46 \\
\num{229376} & 16 & 78 & 82 & 82 & 55 & & 88 & 88 & 54 \\
\num{917504} & 64 & 82 & 78 & 77 & 64 & & 80 & 80 & 63 \\
\num{3670016} & 256 & -- & 75 & 75 & 69 & & 76 & 75 & 69 \\
\bottomrule
\end{tabular}
\label{tab:fp_outer}
\end{table}

The scaling study was repeated using Anderson-acceleration with five Anderson vectors in Table \ref{tab:aa_outer}. 
More uniform convergence with respect to the mesh size is observed when acceleration is employed. 
\begin{table}[h!]
\centering
\caption{Number of outer fixed-point iterations with Anderson acceleration. }
\begin{tabular}{ccccccccccc}
\toprule
 &  &  & \multicolumn{3}{c}{LDG}  &  & \multicolumn{3}{c}{IP} \\
\cmidrule{4-6}\cmidrule{8-10}
$N_e$ & Proc. & P1 & Full & Half & Indep. & & Full & Half & Indep. \\
\midrule
\num{14336} & 1 & 16 & 25 & 26 & 13 & & 25 & 23 & 12 \\
\num{57344} & 4 & 18 & 20 & 21 & 14 & & 22 & 23 & 14 \\
\num{229376} & 16 & 19 & 20 & 20 & 15 & & 20 & 20 & 15 \\
\num{917504} & 64 & 19 & 18 & 18 & 16 & & 18 & 18 & 15 \\
\num{3670016} & 256 & -- & 17 & 17 & 18 & & 17 & 17 & 18 \\
\bottomrule
\end{tabular}
\label{tab:aa_outer}
\end{table}

The scaling of the maximum number of inner, preconditioned conjugate gradient iterations over all outer iterations is shown in Table \ref{tab:weakcp_inner}. 
These results are identical with and without Anderson acceleration. 
Observe that AMG preconditioned conjugate gradient has uniform convergence for all schemes except LDG with half-range boundary conditions. 
Notably the effectiveness of the inner solver was not altered by the choice of boundary conditions for the IP method. 
The half-range boundary conditions introduce a term on the boundary analogous to the term that couples the normal component of the current across interior interfaces in the P$_1$ discretization. 
Recall that this term arises from integrating the gradient of the scalar flux so we surmise that this boundary term introduces a gradient-like scaling in the diagonal block of the LDG discretization that interferes with the effectiveness of AMG. 
The LDG method with full-range boundary conditions and the IP method with either of the boundary conditions are then consistent methods that have scalable solvers comparable to an independent method. 
\begin{table}[h!]
\centering
\caption{Maximum number of inner iterations performed over all outer iterations.}
\begin{tabular}{cccccccccc}
\toprule
 &  & \multicolumn{3}{c}{LDG}  &  & \multicolumn{3}{c}{IP} \\
\cmidrule{3-5}\cmidrule{7-9}
$N_e$ & Proc. & Full & Half & Indep. & & Full & Half & Indep. \\
\midrule
\num{14336} & 1 & 10 & 48 & 10 & & 10 & 11 & 10 \\
\num{57344} & 4 & 11 & 62 & 11 & & 11 & 11 & 11 \\
\num{229376} & 16 & 13 & 86 & 12 & & 12 & 12 & 12 \\
\num{917504} & 64 & 13 & 113 & 13 & & 14 & 13 & 14 \\
\num{3670016} & 256 & 15 & 160 & 15 & & 15 & 16 & 15 \\
\bottomrule
\end{tabular}
\label{tab:weakcp_inner}
\end{table}

\begin{figure}[h!]
\centering
\begin{subfigure}{.32\textwidth}
	\centering
	\includegraphics[width=\textwidth]{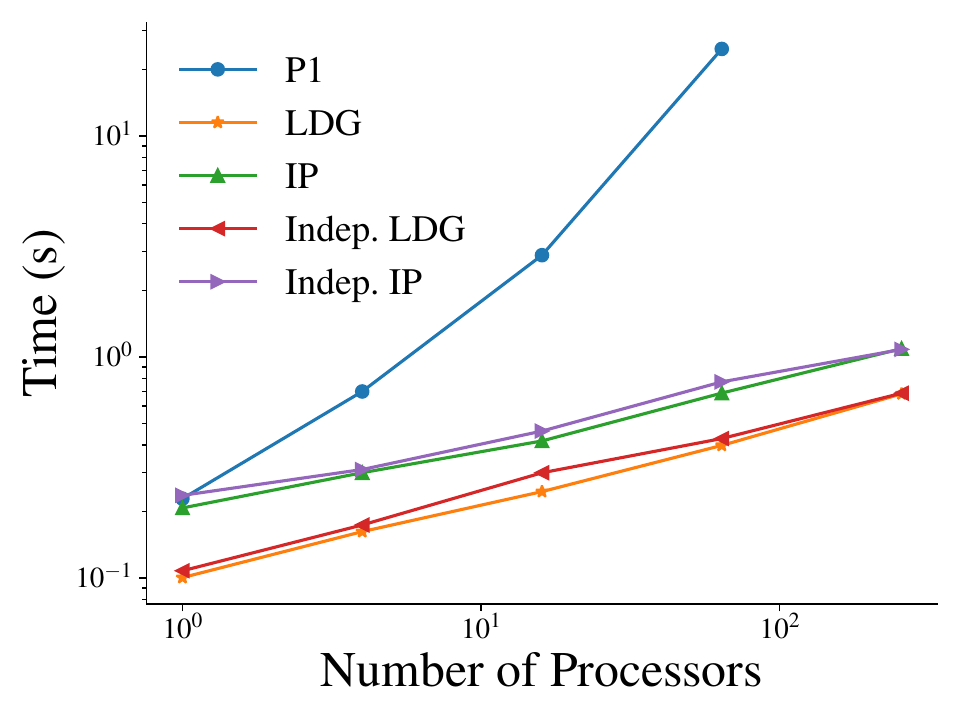}
	\caption{LO solve}
\end{subfigure}
\begin{subfigure}{.32\textwidth}
	\centering
	\includegraphics[width=\textwidth]{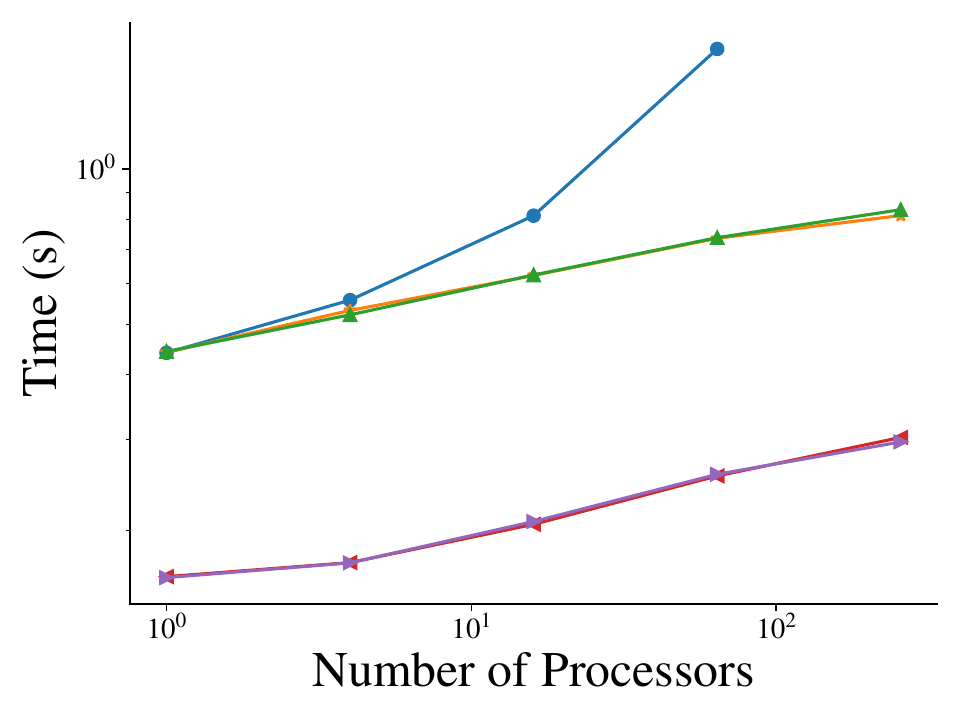}
	\caption{LO assembly}
\end{subfigure}
\begin{subfigure}{.32\textwidth}
	\centering
	\includegraphics[width=\textwidth]{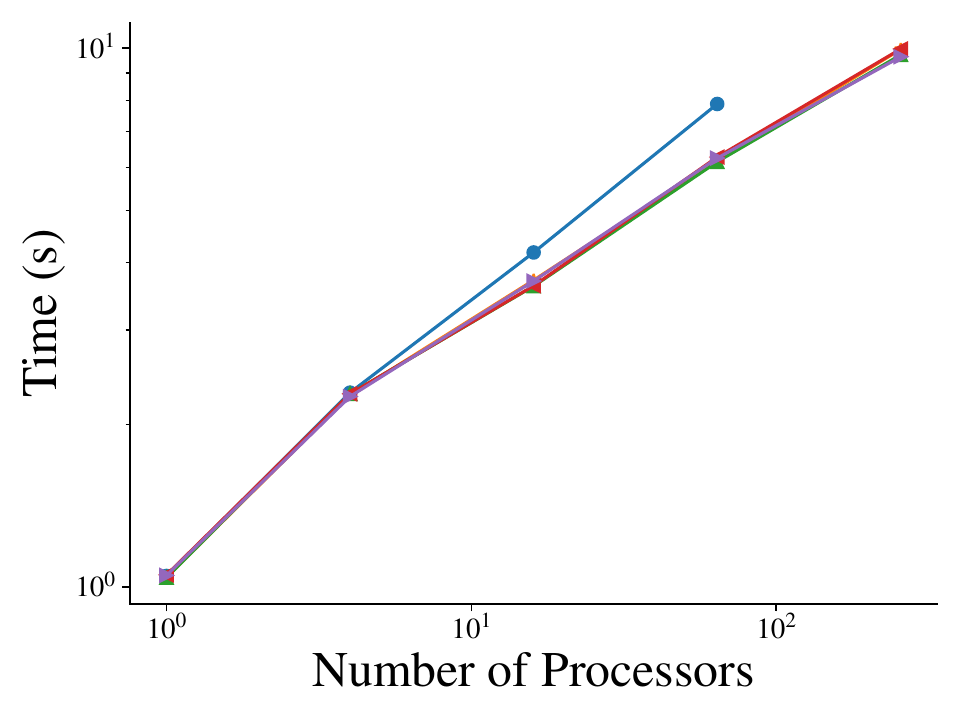}
	\caption{Transport sweep}
\end{subfigure}
\caption{The average cost per iteration to solve the LO system, assemble the LO correction source, and apply the transport sweep. }
\label{fig:cp_cost_per_it}
\end{figure}

The average cost per outer iteration associated with solving the LO system, assembling the SMM correction source, and applying the transport sweep are shown in Fig.~\ref{fig:cp_cost_per_it}. 
Only the consistent LDG and IP methods with full-range boundary conditions are shown since LDG with half-range boundary conditions was not scalable and the cost of IP was insensitive to the choice of boundary condition. 
The sparse direct solver used for the \Pone method is clearly not scalable. 
The fastest LO solutions were obtained by the independent and consistent LDG methods, with the IP methods being about 3\% slower per outer iteration. 
Assembling the consistent SMM source is more expensive than for the independent methods as more computation is needed to form the additional correction source terms required to achieve consistency.
Note that the P$_1$ method shows an unexpectedly poor scaling in the cost of assembly. 
This scaling discrepancy was also seen in the average cost of the transport sweep per iteration. 
The sweep and much of the assembly machinery are shared among all the methods. 
We believe this behavior is due to a possible memory issue with SuperLU as anomalous timings occurred only in a few iterations directly following the LU factorization. 

\begin{figure}[h!]
\centering
\begin{subfigure}{.45\textwidth}
	\centering
	\includegraphics[width=\textwidth]{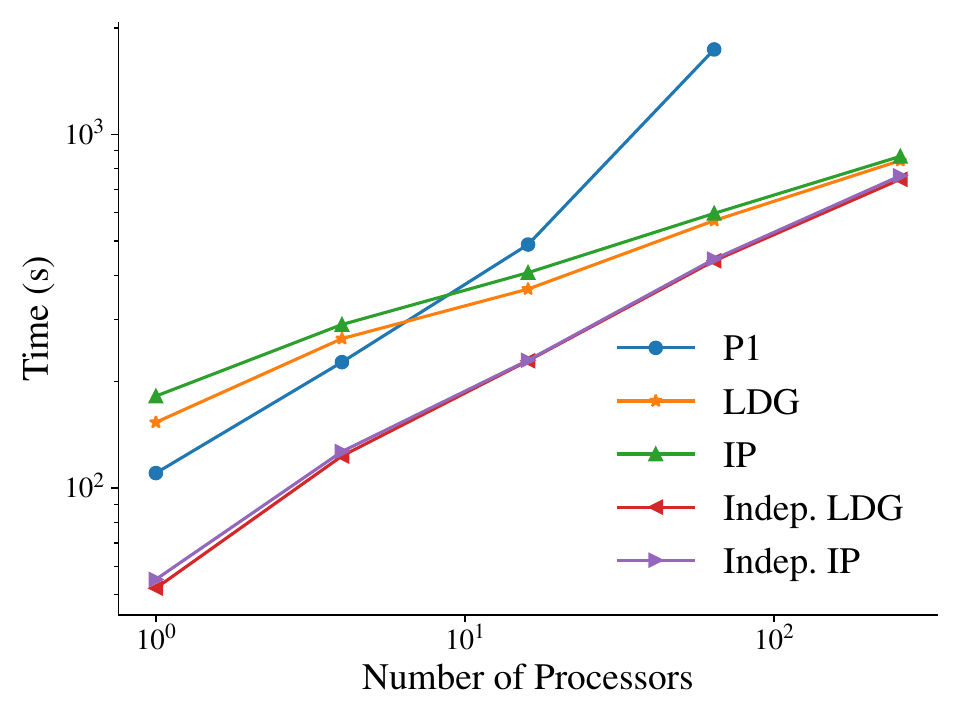}
	\caption{Fixed-point}
\end{subfigure}
\begin{subfigure}{.45\textwidth}
	\centering
	\includegraphics[width=\textwidth]{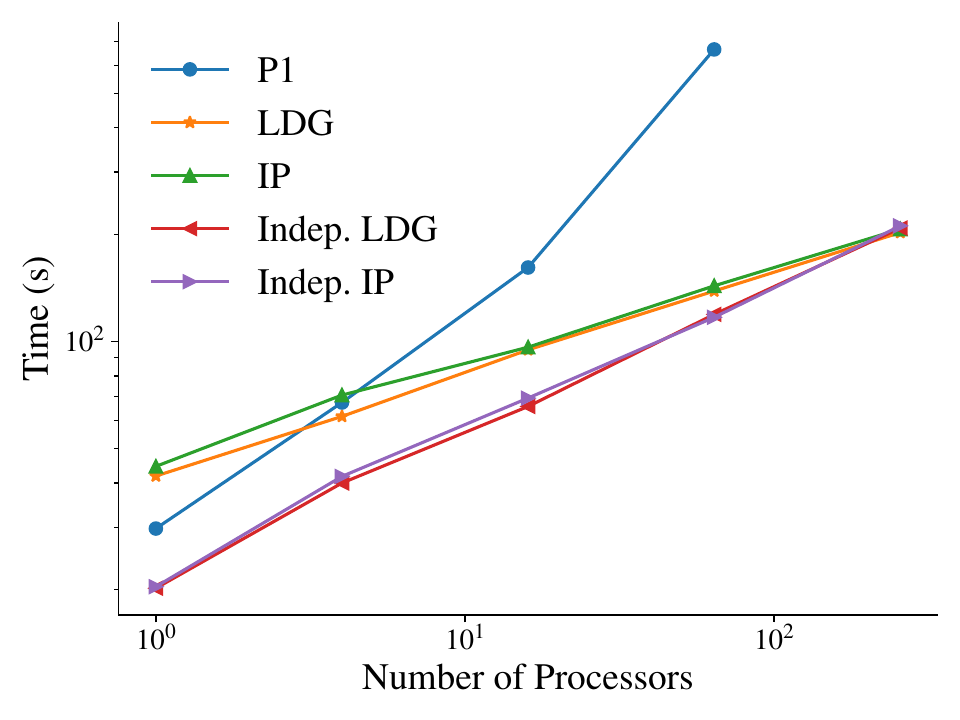}
	\caption{Anderson acceleration}
\end{subfigure}
\caption{Total time to solution for the crooked pipe problem under weak scaling.}
\label{fig:cp_wall}
\end{figure}

The total time to solution is presented in Figure \ref{fig:cp_wall} for solving the crooked pipe with and without Anderson acceleration. 
The effectiveness of the \Pone method in accelerating the outer iteration counterbalances the inefficiency of the LO solver. 
For the first two levels of refinement, the \Pone method is faster than the consistent LDG and IP methods. 
However, the poor parallel scaling of the sparse direct solver makes \Pone significantly more expensive than the consistent LDG and IP methods for subsequent refinements. 
In particular, on the largest mesh for which the direct method could be applied to the \Pone LO system, the consistent LDG and IP methods were up to 3x faster. 
The independent methods were the fastest at each mesh resolution, resulting in speedups of 13\% over the scalable, consistent variants on the finest mesh. 

With Anderson acceleration, the cost of factorizing the P$_1$ system is amortized across fewer outer iterations. 
This results in up to a 5x speedup of consistent LDG and IP over P$_1$. 
Since Anderson narrows the gap in outer iterations between consistent and independent, consistent was able to beat independent on the most resolved problem. 
Thus, when the discretization error is small enough, the proposed consistent methods produce a consistent solution while achieving comparable performance to independent methods. 
Nuances in cost of forming and solving the LO system are dominated by the high cost of the transport sweep, resulting in all the scalable algorithms producing similar runtimes on the most refined mesh. 

Using the data from the Anderson-accelerated solve, the LDG and IP-based methods exhibit weak scaling efficiency of the inner, linear solve between 10\% and 20\% on 256 processors. 
The IP-based methods scaling was slightly better than the LDG-based methods. 
The \Pone moment solve did not scale due to the use of the sparse direct method. 
The full algorithms scale between 15\% and 30\% for the IP and LDG-based methods with \Pone resulting in an overall unscalable transport algorithm. 
Note that the weak scaling of the IP and LDG diffusion systems is less than what was reported in \citet{osti_1117924} for \texttt{BooomerAMG} applied to similar problems. 
This is likely because their weak scaling is compared to a base case starting with 64 processors while our base case is serial. 
Using a serial base case makes the inclusion of parallel communication costs in subsequent calculations more pronounced. 
We did observe similar iterative convergence comparable to what was seen in \cite{osti_1117924}.

\subsubsection{Comparison of Solution Quality}
\begin{figure}[h!]
\centering
\begin{subfigure}{.49\textwidth}
	\centering
	\includegraphics[width=\textwidth]{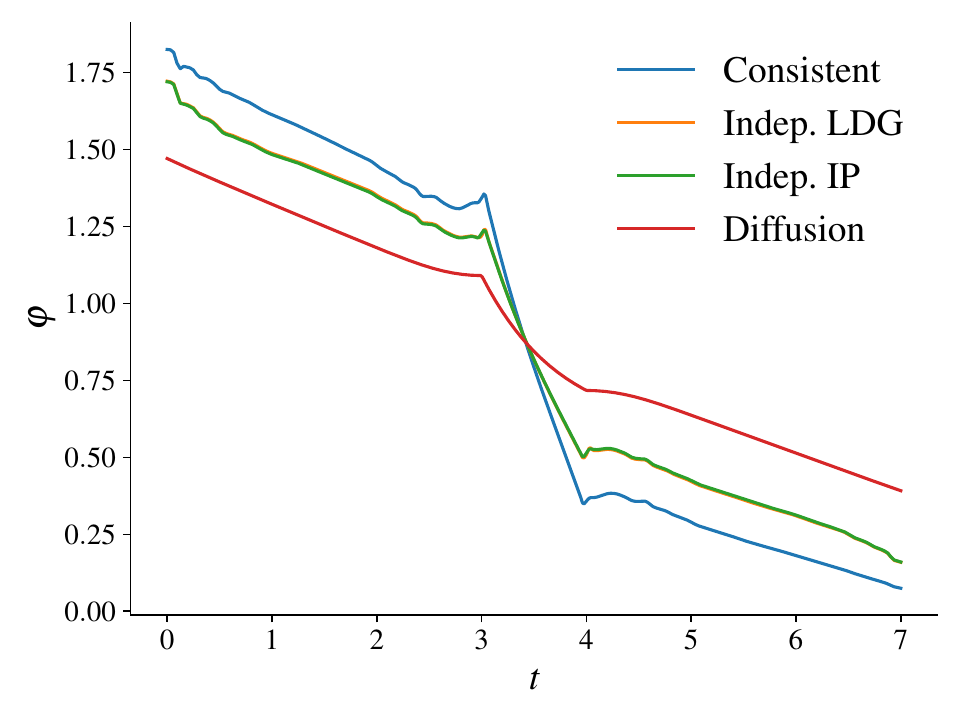}
	\caption{$N_e = \num{14336}$}
\end{subfigure}
\begin{subfigure}{.49\textwidth}
	\centering
	\includegraphics[width=\textwidth]{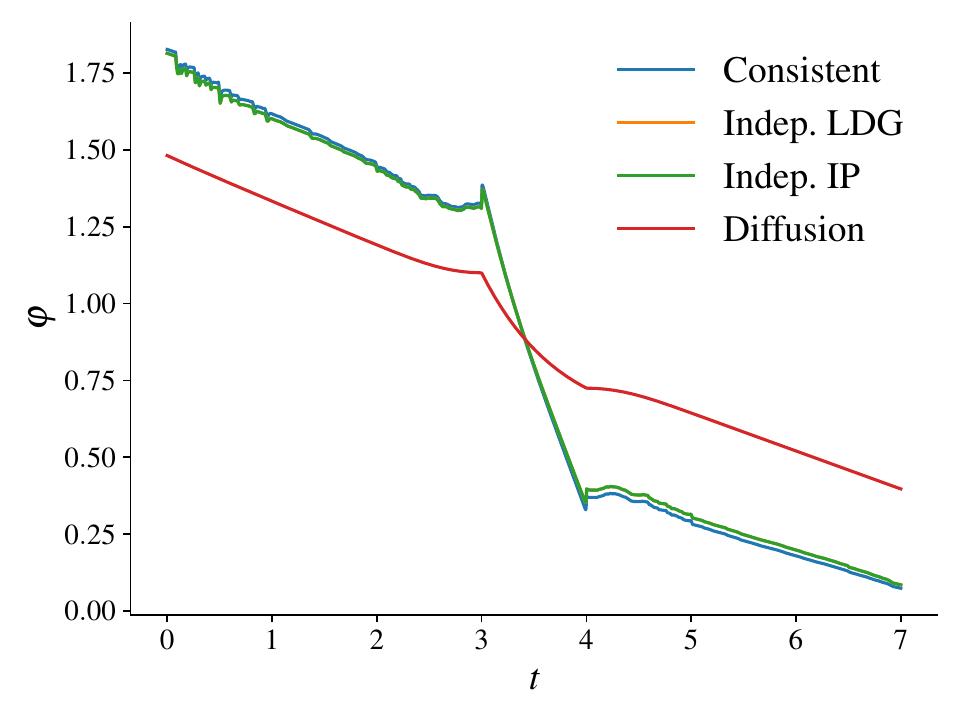}
	\caption{$N_e = \num{3670016}$}
\end{subfigure}
\caption{Scalar flux along $y=0$. }
\label{fig:line_phi}
\end{figure}

Solutions generated by the consistent methods are compared with those of the independent methods on the base mesh and on the most refined mesh. 
The consistent methods attain the same solutions up to the $10^{-6}$ fixed-point iteration tolerance so only the consistent LDG solution is shown. 
Figure \ref{fig:line_phi} shows a lineout of the scalar flux solution along the centerline $y=0$ on the least and most refined meshes. 
A diffusion solution on the same mesh provides a reference with which to compare the physics fidelity of the SMMs. 
On the least refined mesh, the independent methods have numerical diffusion that results in under and over heating the front and back of the inner wall, respectively, as compared to the consistent solution. 
This numerical diffusion pushes the independent methods closer toward the diffusion solution. 
On the refined mesh, independent and consistent produce similar solutions. 

\begin{figure}[h!]
\centering
\begin{subfigure}{.49\textwidth}
	\centering
	\includegraphics[width=\textwidth]{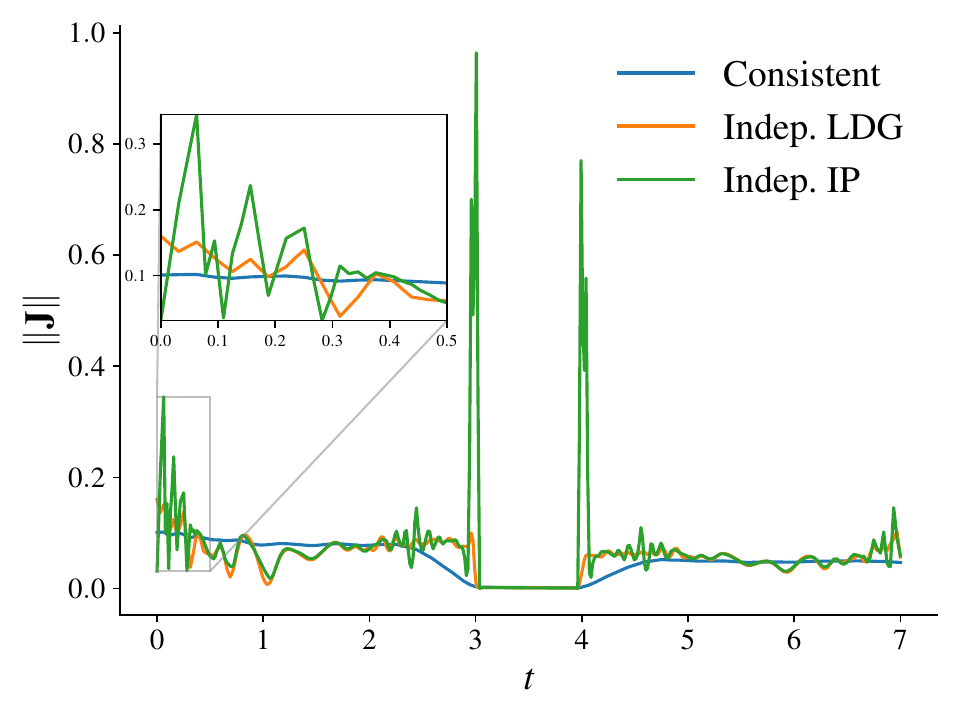}
	\caption{$N_e = \num{14336}$}
\end{subfigure}
\begin{subfigure}{.49\textwidth}
	\centering
	\includegraphics[width=\textwidth]{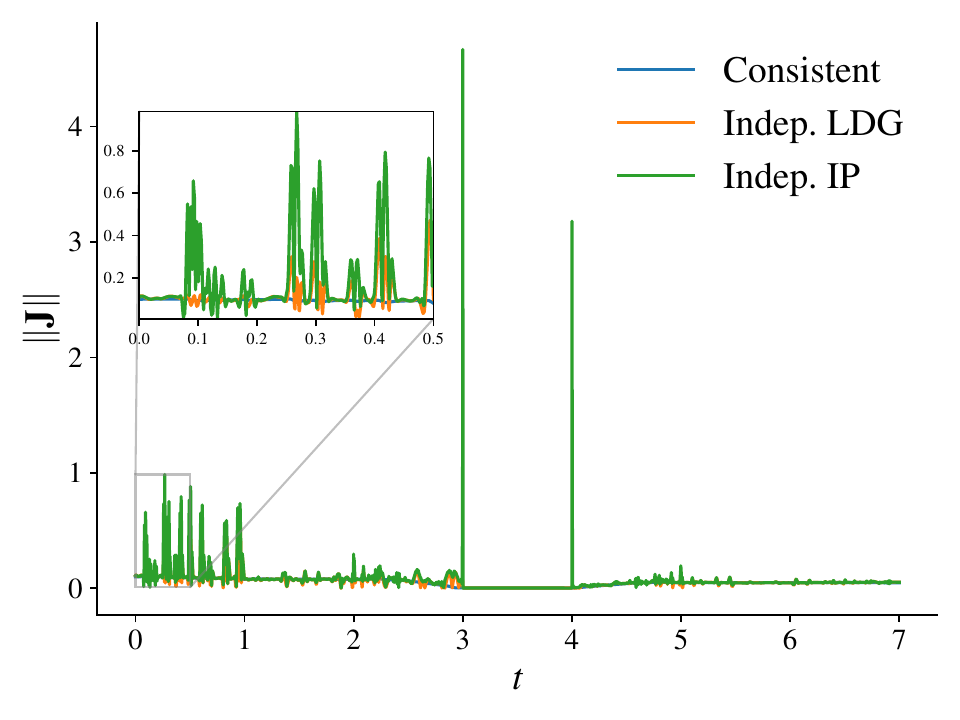}
	\caption{$N_e = \num{3670016}$}
\end{subfigure}
\caption{Magnitude of the current along the line $y=0$.}
\label{fig:line_J}
\end{figure}

Plots of the magnitude of the current along the same centerline $y=0$ are shown in Fig.~\ref{fig:line_J}. 
The independent methods produce currents with wide, non-physical oscillations. 
These oscillations increase in magnitude as the mesh is refined, particularly at the wall-pipe interface. 
By contrast, the consistent scheme achieves a physically realistic solution on both meshes. 
It is possible that better discretizations of the LO system, such as Raviart Thomas, could improve the solution quality of an independent SMM. 
However, a Raviart Thomas-based SMM was also shown to have suboptimal accuracy of the current in \cite{olivier_smm}. 

\begin{figure}[h!]
\centering
\begin{subfigure}{.35\textwidth}
	\centering
	\includegraphics[width=\textwidth]{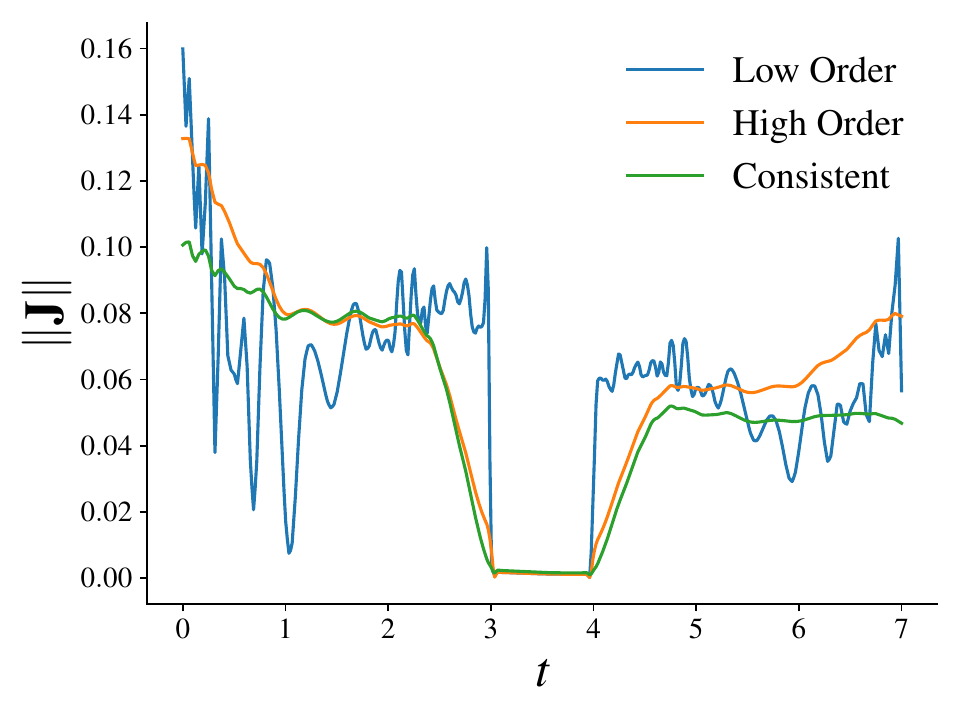}
	\caption{LDG $N_e = \num{14336}$}
\end{subfigure}
\begin{subfigure}{.35\textwidth}
	\centering
	\includegraphics[width=\textwidth]{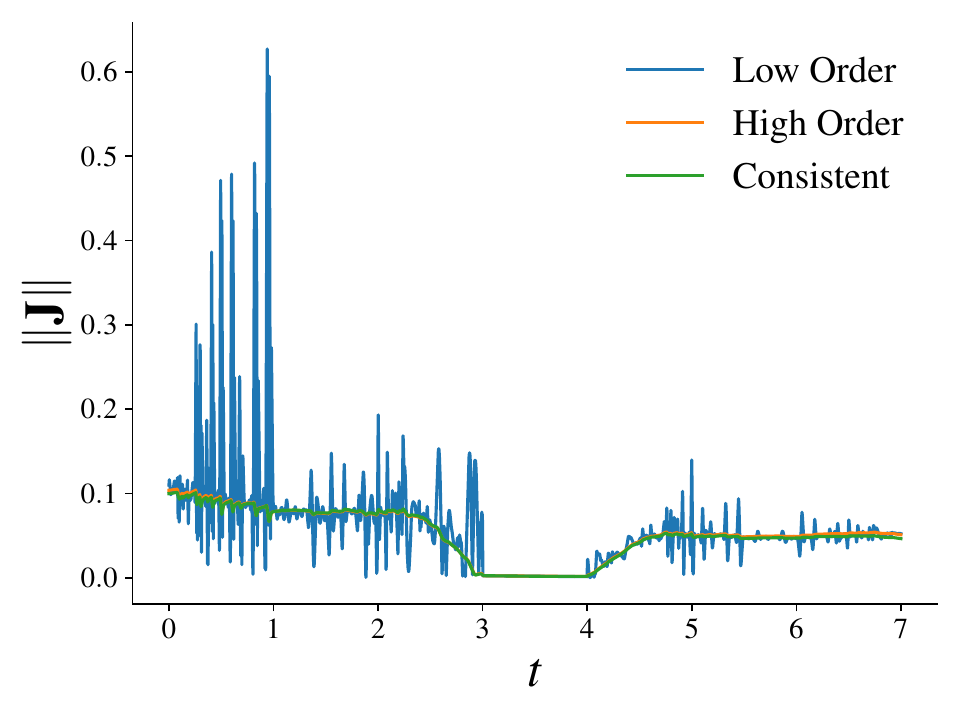}
	\caption{LDG $N_e = \num{3670016}$}
\end{subfigure}
\begin{subfigure}{.35\textwidth}
	\centering
	\includegraphics[width=\textwidth]{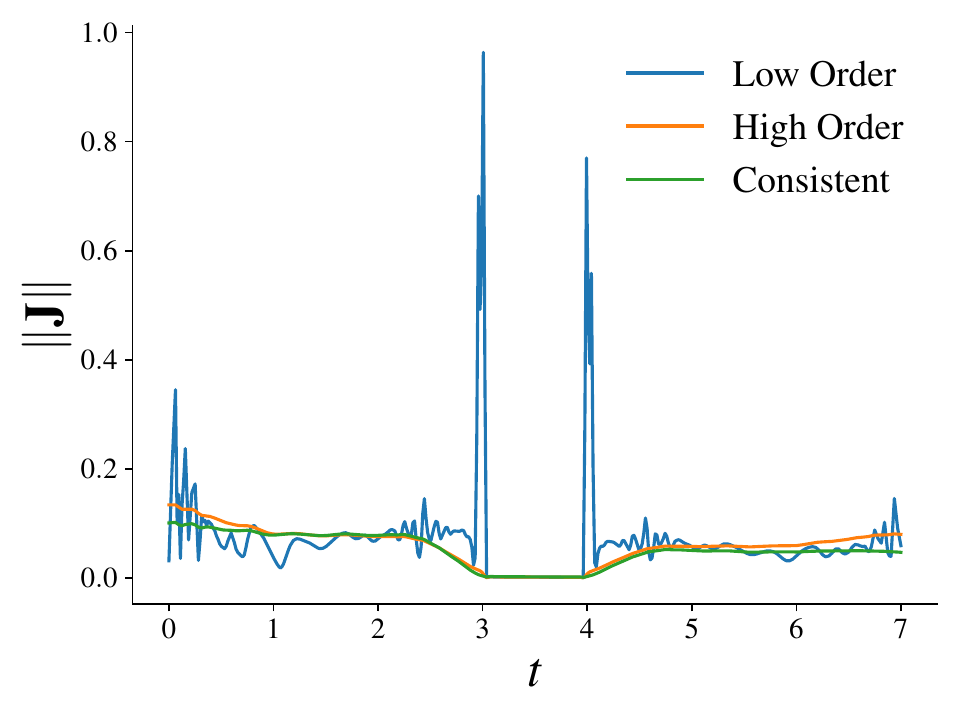}
	\caption{IP $N_e = \num{14336}$}
\end{subfigure}
\begin{subfigure}{.35\textwidth}
	\centering
	\includegraphics[width=\textwidth]{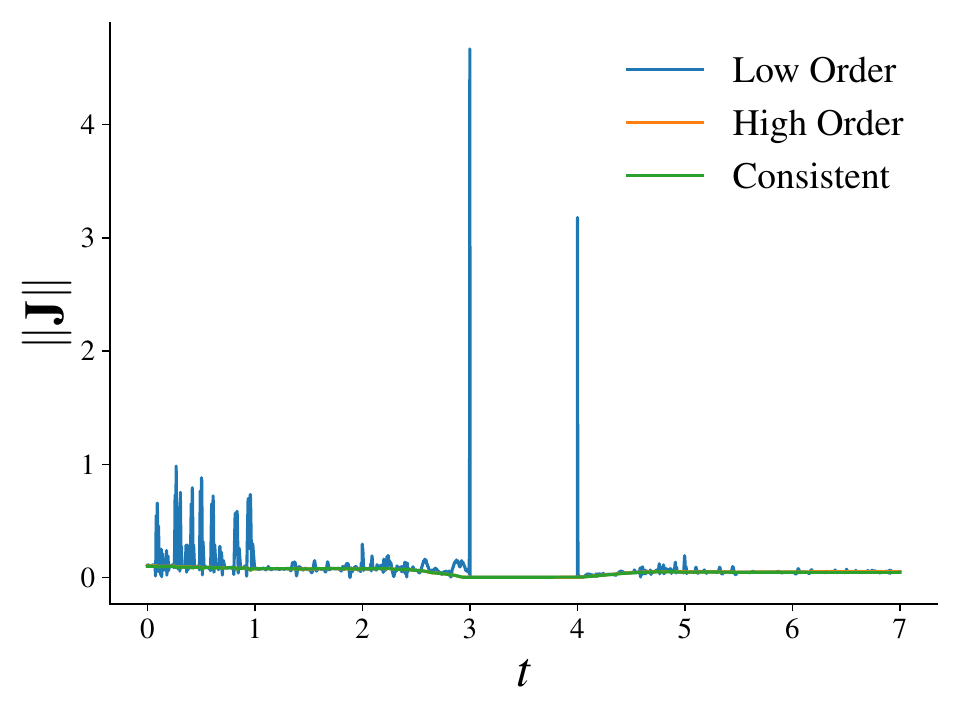}
	\caption{IP $N_e = \num{3670016}$}
\end{subfigure}
\caption{Magnitude of the HO and LO current for the independent SMM along $y=0$.}
\label{fig:indep_ho_v_lo}
\end{figure}

In Fig.~\ref{fig:indep_ho_v_lo}, the HO and LO currents from the independent methods are compared to the consistent current solution. 
The HO current associated with the independent methods does not exhibit the wild oscillations that the independent LO currents do on both meshes. 
However, use of the HO solution may have implications for conservation, especially in multiphysics simulations. 
Consistent methods avoid such issues by producing HO and LO solutions that are both conservative. 

\section{Conclusions}\label{sec:conc}
We developed three novel consistent Second Moment Methods (SMMs) based on low-order (LO) systems with left hand sides equivalent to the fully consistent \Pone, Local Discontinuous Galerkin (LDG), and Interior Penalty (IP) discretizations of radiation diffusion. 
A discrete residual approach was used to derive the SMM correction source terms that make these choices of LO diffusion discretizations consistent with the Discontinuous Galerkin (DG) discretization of the Discrete Ordinates (\Sn) transport equations. 
The consistent SMMs were compared to the independent LDG and IP SMM variants. 
We found that the consistent SMMs were more accurate than the independent methods. 
In particular, the currents generated by the consistent methods converged with optimal second-order accuracy while the independent LDG and IP methods showed $\mathcal{O}(h^{3/2})$ and $\mathcal{O}(h)$ accuracy, respectively. 
The LDG and IP methods were scalably solved in parallel with algebraic multigrid, obtaining a LO solution consistent with the DG \Sn transport discretization on a challenging, multi-material benchmark problem. 
The \Pone method was effective at accelerating the outer fixed-point iteration but the LDG and IP methods were more efficient on more resolved problems because the sparse direct method used to solve the \Pone LO system was not scalable. 
Generally, the independent variants of LDG and IP were the most efficient methods, requiring fewer total transport sweeps than the consistent variants. 
This discrepancy was less apparent as the mesh was refined because the additional discretization-dependent correction terms for the consistent methods become smaller in magnitude. 
The efficiency of the independent methods came at the cost of poor solution quality for the LO solution, especially for the current. 
We therefore recommend the use of the LDG or IP consistent SMMs over their independent variants. 

\section{Acknowledgements}
This research used resources provided by the Darwin testbed at Los Alamos National Laboratory (LANL) which is funded by the Computational Systems and Software Environments subprogram of LANL's Advanced Simulation and Computing program. 
LANL is operated by Triad National Security, LLC, for the National Nuclear Security Administration of the U.S. Department of Energy (Contract No.~89233218CNA000001). 

\bibliographystyle{IEEEtranN}
\bibliography{references}

\end{document}